\documentclass[a4paper,11pt]{amsart}

\usepackage{amsfonts}
\usepackage{latexsym}
 \usepackage{amssymb}
\usepackage[all]{xy}

\renewcommand{\Sigma}{\sum}
\newtheorem{thm}{\sc theorem}[section]
\newtheorem{dfn}[thm]{\sc definition}

\newtheorem{lemma}[thm]{\sc lemma}

\newtheorem{cor}[thm]{\sc corollary}

\newenvironment{pf}{{{\bf Proof.}}}{\hfill$\Box$\\[1mm]}

\newcommand{\A}{{\mathcal{A}}}
\newcommand{\dA}{{\mathsf{d}\mathcal{A}}}
\newcommand{\B}{{\mathcal{B}}}
\newcommand{\BB}{{\mathsf{B}}}

\renewcommand{\b}{\mathsf{b}}

\renewcommand{\k}{{\mathbb{C}}}

\newcommand{\M}{{\mathcal{M}}}
\newcommand{\N}{{\mathcal{N}}}
\newcommand{\NN}{{\mathsf{N}}}
\renewcommand{\r}{\mathbf{r}} 
\newcommand{\s}{\mathsf{s}} 
\renewcommand{\t}{\mathsf{t}} 
\newcommand{\T}{\mathsf{T}} 
\newcommand{\U}{{\mathcal{U}}}

\newcommand{\Slambda}{S(\lambda)}

\newcommand{\ccc}{\xi} 
\newcommand{\half}{{\frac{1}{2}}}

\newcommand{\even}{\mathrm{even}}
\newcommand{\odd}{\mathrm{odd}}
\newcommand{\otherwise}{\mathrm{otherwise}}

\begin{document}

\title{Twisted Homology of Quantum $SL(2)$ - Part II}

\author{Tom~Hadfield} 
\author{Ulrich Kr\"ahmer} 

\address{Gloucester Research, Whittington House, 19-30 Alfred Place, London WC1E 7EA, UK}
\email{Thomas.Daniel.Hadfield@gmail.com}

\address{University of Glasgow, Department of Mathematics\\University Gardens, Glasgow G12 8QW, UK}
\email{ukraehmer@maths.gla.ac.uk}

\begin{abstract}
We
complete the calculation of the twisted cyclic 
homology of the quantised
coordinate ring $\A=\k_q[SL(2)]$ of $SL(2)$
that we began in \cite{kt}. In particular, a 
nontrivial cyclic 3-cocycle is constructed which also has
a nontrivial class in Hochschild cohomology and thus should be
viewed as a noncommutative geometry analogue of a volume form. 
\end{abstract}


\maketitle

\centerline{\emph{Dedicated to Prof.~K.~Schm\"udgen on the
occasion of his 60th birthday}}
 
\section{Introduction}
In this paper we finish the
computation \cite{kt} of the twisted cyclic homology 
$HC_\bullet^\sigma(\A)$ of 
the quantised coordinate ring 
$\A=\k_q[SL(2)]$ for
$q \in \k$ not a root of unity. 
This invariant was defined by
Kustermans, Murphy and Tuset
\cite{kmt} and modifies Connes' cyclic
homology
$HC_\bullet(\A)=HC_\bullet^{\mathrm{id}}(\A)$
using an automorphism 
$\sigma$ of $\A$. 
It is computed using the Connes spectral sequence 
\begin{equation}\label{css}
		  E^1_{pq}=H_{q-p}(\A,{}_\sigma \A) \Rightarrow 
		  HC_{p+q}^\sigma (\A),
\end{equation}
where $H_\bullet(\A,\,\cdot\,)$
is the Hochschild homology of $\A$
and the coefficient bimodule ${}_\sigma \A$ arises from $\A$
by twisting the canonical left action to 
$a \triangleright b=\sigma(a)b$.

In \cite{kt} we were unable to fully compute
the spectral sequence in the ``no dimension drop''
case, that is, for those $\sigma$ with 
$H_3(\A,{}_\sigma \A) \neq 0$.
 The crucial
new ingredient we use in this paper are the cup and cap products
\begin{eqnarray}
&& \smallsmile : H^m(\A,{}_\sigma \A) \otimes 
		  H^n(\A,{}_\tau \A) \rightarrow 
		  H^{m+n}(\A,{}_{\tau \circ \sigma}
		  \A),\nonumber\\ 
&& \smallfrown : H_n(\A,{}_\sigma \A) 
		  \otimes H^m(\A,{}_\tau \A) \rightarrow
		  H_{n-m} (\A,{}_{\tau^{-1} \circ \sigma}
		  \A).\nonumber
\end{eqnarray} 
together with the twisted derivations
$$
		  \partial^\pm_H,\partial^\pm_E,\partial^\pm_F
		  \in \bigoplus_{\sigma \in
		  \mathrm{Aut}(\A)} H^\bullet(\A,{}_\sigma
		  \A)
$$
that deform the action of left- and
right-invariant vector fields on $SL(2)$. 
In Section~3 we show that under $\smallsmile$  these derivations generate a 
$q$-deformed exterior algebra whose
$\smallfrown$-action allows us to
identify nontrivial 
2- and 3-cycles. Using this we then compute
$HC_\bullet^\sigma (\A)$ in Section~\ref{schnoel}.

Our motivation 
for \cite{kt} was the relation of 
$HC_\bullet^\sigma(\A)$
to Woronowicz's
theory of covariant differential calculi \cite{kmt}.
We realised  subsequently 
\cite{cr} that the
coefficients ${}_\sigma \A$ also arise from
Poincar{\'e} duality in Hochschild
(co)homology: using the general theory
of Van den Bergh \cite{vdb} we showed
that
\begin{equation}\label{edin}
		  H^n(\A,{}_\sigma \A) \simeq
			H_{3-n}(\A,{}_{\sigma_{q^{-2},1}
			\circ \sigma} \A) \quad\forall\,
			\sigma \in \mathrm{Aut} (\A),
\end{equation} 
where for $ \lambda , \mu \in \k \setminus \{0\}$, 
$ \sigma_{\lambda,\mu} \in \mathrm{Aut} (\A)$ is
determined by its values 
$$ \sigma_{\lambda,\mu} (a)=\lambda a,\quad
\sigma_{\lambda,\mu} (b)=\mu b,\quad
\sigma_{\lambda,\mu}(c)=\mu^{-1} c,\quad
\sigma_{\lambda,\mu}(d)=\lambda^{-1}d
$$ 
on the standard
generators $a,b,c,d$ of $\A$. In particular,
$\sigma_{q^{-2},1}$ is Woronowicz's modular automorphism of the
Haar functional on $\A$, and the ``no
dimension drop'' automorphisms are
precisely 
$\sigma = \sigma_{q^{-N},1}$,
 with $N \in \mathbb{Z}$, 
 $N \geq 2$.

Our main results are summarised in
Theorems~\ref{themainresult} and
\ref{themainresult2} below. Therein,
$$\dA \in H_3(\A,{}_{\sigma_{q^{-2},1}} \A)$$ 
is the fundamental class in Hochschild
homology, that is,  it corresponds under (\ref{edin})
to $1 \in H^0(\A,\A)$ (identified with the centre of
$\A$), and (\ref{edin})
is given by $\cdot \smallfrown \dA$ \cite{nielsuli}. 
A cycle in the standard Hochschild 
complex $C_\bullet(\A,{}_{\sigma_{q^{-2},1}} \A)=
\A^{\otimes \bullet+1}$ representing
$\dA$ is given explicitly in 
 (\ref{defda}).

As a consequence of (\ref{edin}) 
the Connes spectral sequence (\ref{css}) stabilises at the second page, 
and $HC_n^\sigma (\A) \simeq HC_{n+2}^\sigma (\A)$
for $n \ge 3$. 
In this way, we obtain the two periodic
cyclic homology groups 
 $HP^{\sigma}_\even(\A)$
and $HP^{\sigma}_\odd(\A)$
as the limit of the
$HC_{2n}^\sigma (\A)$ and
$HC_{2n+1}^\sigma (\A)$, respectively. 
These groups can be described 
in the ``no dimension drop'' case left
open in \cite{kt} as follows:
\begin{thm}\label{themainresult}
Let $\A$ be the quantised coordinate ring $\k_q[SL(2)]$,
$q \in \k$ not a root of unity, and 
$ \sigma=\sigma_{q^{-N},1} \in
 \mathrm{Aut}(\A)$, $N \in \mathbb{Z}$,
be determined by
$$
		  \sigma(a)=q^{-N}a,\quad
		  \sigma(b) = b.			  
$$
Then the $ \sigma$-twisted periodic cyclic
homology of $\A$ is given by
\begin{eqnarray}
&& \!\!\!\!\!\!\!\!\!\! HP^\sigma_\even(\A)=
		  {HC^\sigma_4}(\A)=\left\{
		  \begin{array}{ll}
		  \k[b^rc^r (\dA \smallfrown [{\partial^-_H}])] 
			& : N=2r+2, r \ge 0,\\
		  \k[1] & : \otherwise,
		  \end{array}\right.  \nonumber\\ 
&&  \!\!\!\!\!\!\!\!\!\! HP^\sigma_\odd (\A)=
		  HC^\sigma_3 (\A)=\left\{
		  \begin{array}{ll}
		  \k[b^rc^r \dA] 
			\quad & :\,N=2r+2,\,r \ge 0,\\
		  \k[b \otimes c] \quad & :\,\otherwise.
		  \end{array}\right.\nonumber  
\end{eqnarray} 
Here classes in $HC_n^\sigma (\A)$ 
are represented by
 classes in 
$H_{n-2p}(\A,{}_\sigma\A)$, $p \ge 0$ using Connes'
spectral sequence 
$E^1_{pq}=H_{q-p}(\A,{}_\sigma\A) \Rightarrow 
HC_{p+q}^\sigma (\A)$.
\end{thm}

Our second main result, which we prove in  Sections ~\ref{como} and \ref{cyccohom},
  is the explicit construction of a twisted cyclic 3-cocycle $\ccc$
 that pairs nontrivially
with $\dA$:
\begin{thm}\label{themainresult2}
Define for all $j,k \ge 0$ a functional 
$\int_{[b^jc^k]} : \A \rightarrow \k$ by
$$
		  \int_{[b^jc^k]} e_{r,s,t}:=\delta_{0,r} 
		\delta_{j,s} \delta_{k,t},\quad
		  e_{i,j,k}:=\left\{
		  \begin{array}{ll}
		  a^i b^j c^k\quad & :\;i \ge 0\\
		  d^{-i}b^jc^k\quad & :\;i<0
		  \end{array}\right.
$$ 
and two linear functionals 
$C_3(\A,{}_{\sigma_{q^{-2},1}} \A) \rightarrow \k$ by
\begin{eqnarray}
&& \varphi(\cdot):=
			\int_{[1]} 
			\cdot \smallfrown (\partial^+_H \smallsmile
			\partial^+_E \smallsmile
			\partial^+_F),\nonumber\\ 
&& \eta(\cdot):=2\int_{[bc]}
		  \cdot \smallfrown
		  (\partial^+_H \smallsmile 
		  (\sigma_{1,1/2}-\mathrm{id})\smallsmile 
		  (\sigma_{1,2}-\mathrm{id})).\nonumber
\end{eqnarray} 
Then $\varphi$ and $\eta$ are 
$ \sigma_{q^{-2},1}$-twisted Hochschild 3-cocycles,
$\ccc := \varphi +\eta$ is a 
$ \sigma_{q^{-2},1}$-twisted cyclic
3-cocycle, 
 and
$\ccc(\dA)= 1 = \varphi(\dA)$.
\end{thm}

The 3-cocycle $\ccc$ 
gives an explicit description in terms of 
Connes' $\lambda $-complex of 
$HC^3_{\sigma_{q^{-2},1}}(\A)=
(HC_3^{\sigma_{q^{-2},1}}(\A))^* \simeq 
\k$ and should be interpreted from the
point of view of noncommutative geometry
as an analogue of a volume form, compare 
e.g.~\cite{connes1}, Corollary~35 on
p.337 and Connes' orientability axiom
\cite{cbuchn}, Definition~1.233 
which requires that the Chern character
of a spectral triple should have a
nontrivial class in Hochschild cohomology.
It had been shown by Masuda, Nakagami
and Watanabe that there is no such
volume form in the untwisted $HC^3(\A)$
\cite{mnw} but only in $HC^1(\A)$, and
this ``dimension drop'' was considered
by many authors to be ``rather esoteric''
\cite{connes7}. The cocycle $\xi$ solves
that mystery, so we expect that
there is a spectral triple that realises
$\xi$ via a twisted variant of the Connes-Moscovici
local index formula, for example as in \cite{triest,nestus,sw}.

As shown by Etingof and
Dolgushev \cite{etingof}, 
 a Poincar\'e-type duality as in (\ref{edin})
holds for formal deformation quantisations of
smooth Poisson varieties, and also, as shown by Brown
and Zhang, for a large class of Noetherian
Hopf algebras that includes in particular the quantised
coordinate rings $\k_q[G]$
for all simple algebraic groups $G$ 
\cite{bz}. See also \cite{farinati,uli}
for more information.

Most of the computations in this paper have been verified with
the help of the computer algebra system 
FELIX \cite{felix}.
The FELIX output is
available in electronic form \cite{code}, and it
can easily be adapted to perform similar
computations with other algebras given in terms of
generators and relations.

\subsection*{Acknowledgements}

T.H. and U.K.: We thank S.~Launois for
kindly pointing out to us 
an inconsistency
between our results in \cite{kt} and \cite{cr}, 
which arises from a mistake in our computation 
of $H_2(\A,{}_\sigma \A)$ in \cite{kt} (see \cite{LL}). 
We also thank the referees for their careful reading of the paper and their many helpful and insightful comments.

U.K.: My work was supported by the Marie Curie
fellowship EIF 515144 and the EPSRC 
fellowship EP/E/043267/1. 
Further thanks go to Istv\'an Heckenberger 
who introduced me to the computer algebra 
system FELIX, and to Andreas Thom, Boris Tsygan,
and Christian Voigt for discussions.  

\section{Hochschild homology}\label{hoho}
\subsection{Background}\label{hoho1} 
In Section~\ref{hoho1}, we fix notation and
conventions concerning 
Hochschild homology and the quantised
coordinate ring of $SL(2)$.
For details and proofs see for example 
\cite{CE,loday,weibel} and \cite{joseph,chef} respectively.

\subsubsection{Algebras and bimodules} 
Throughout this paper 
``algebra'' means unital associative $\k$-algebra. 
An unadorned $ \otimes $ denotes the tensor product of
$\k$-vector spaces. 
For an $\A$-bimodule $\M$ and two automorphisms 
$\sigma , \tau $ of an algebra $\A$ we denote by
${}_\sigma \M_\tau $ the bimodule 
which is $\M$ as vector space with bimodule
structure $ x \blacktriangleright 
y \blacktriangleleft z:=
\sigma (x) \triangleright y \triangleleft 
\tau (z)$, $x,z \in \A,y \in \M$, where 
$ \triangleleft , \triangleright $ are the original
actions on $\M$. Note the bimodule isomorphisms
$$
{}_{\sigma'}({}_\sigma \M_\tau)_{\tau'} \simeq 
{}_{\sigma \circ \sigma'} \M_{\tau \circ \tau'},\quad
\M \otimes_\A {}_\sigma \N \simeq
\M_{\sigma^{-1}} \otimes_\A \N,\quad
\A_{\sigma^{-1}} \simeq {}_\sigma \A.
$$

\subsubsection{Hochschild homology}\label{china}  
The Hochschild homology groups 
of an algebra 
$\A$ with coefficients in an $\A$-bimodule 
$\M$
are 
$$
		  H_n(\A,\M) =\mathrm{Tor}_n^{\A^e}(\M,\A),
$$
where $\A^e:=\A \otimes \A^\mathrm{op}$ is the
enveloping algebra of $\A$ (so $\A^e$-modules are just
$\A$-bimodules). They can be computed using the canonical 
bar resolution of $\A$ as an $\A^e$-module, 
and are then realised as the simplicial homology of
the simplicial $\k$-vector space
$C_\bullet(\A,\M):=\M \otimes \A^{\otimes \bullet}$
whose structure maps are
\begin{eqnarray}
&& \b_0 : a_0 \otimes \ldots \otimes a_n \mapsto
	a_0 \triangleleft a_1 \otimes a_2 \otimes 
	\ldots \otimes a_n,\nonumber\\ 
&& \b_i : a_0 \otimes \ldots \otimes a_n \mapsto
	a_0 \otimes \ldots \otimes a_ia_{i+1} \otimes
	\ldots \otimes a_n,\quad 0<i<n,\nonumber\\ 
&& \b_n : a_0 \otimes \ldots \otimes a_n \mapsto
	a_n \triangleright a_0 \otimes \ldots \otimes a_{n-1},
	\nonumber\\ 
&& \s_i : a_0 \otimes \ldots \otimes a_n \mapsto
	a_0 \otimes \ldots \otimes a_i \otimes 1 \otimes 
	a_{i+1} \otimes \ldots \otimes a_n,\quad 0 \leq i \leq n.\nonumber  
\end{eqnarray} 
That is $H_\bullet(\A,\M)$ is (isomorphic to)
the homology of the chain complex 
$C_\bullet(\A,\M)$ 
whose boundary map is 
given by
$$
		  \b := \sum_{i=0}^{n} (-1)^i \b_i.
$$
In the sequel we will often write  
$\b(a_0,\ldots,a_n)$ instead of
$\b(a_0 \otimes \ldots \otimes a_n)$.

If $\A$ is the coordinate ring $\k[X]$ of a smooth affine
variety, then $H_\bullet(\A,\A)$ can be identified
canonically with the K\"ahler differentials (algebraic
differential forms) on $X$
(Hochschild-Kostant-Rosenberg theorem). 

\subsubsection{The normalised complex}
For any simplicial $\k$-vector space $C$, 
$$
		  D:=\mathrm{span}\{\mathrm{im}\,\s_i\} \subset C
$$
is a contractible subcomplex with respect to $\b$, 
so the canonical map
$$
		  C \rightarrow \bar C:=C/D
$$ 
is a quasi-isomorphism of complexes, and  
 working with the so-called 
normalised complex $(\bar C_\bullet,\b)$ 
simplifies many computations.

For 
$C=C(\A,\M)$ from the previous section, 
we have $\bar C_n = \M \otimes \bar \A^{\otimes n}$,
where $\bar \A:=\A/\k$. 
So informally speaking in the
computation of Hochschild homology we can neglect
 all elementary
tensors with a tensor component being equal to a
multiple of $1 \in \A$. 

\subsubsection{Quantum $SL(2)$}\label{aldlos} 
For the remainder of Section~\ref{hoho}, 
$q \in \k$ denotes a fixed nonzero
parameter, which we assume is not a root of unity. 
Furthermore, $\A$ is throughout the quantised 
coordinate algebra 
$\k_q[SL(2)]$ of $SL(2)$, that is the algebra generated by
symbols $a$, $b$, $c$, $d$ with relations
$$
		  ab=qba, \quad
		  ac=qca, \quad
		  bc=cb,\quad
		  bd=qdb, \quad
		  cd=qdc,
$$
$$
		  ad-qbc=1, \quad
		  da - q^{-1} bc=1.
$$
It follows from the defining relations that the elements
\begin{equation}\label{eijk}
		  e_{i,j,k}:=\left\{
		  \begin{array}{ll}
		  a^i b^j c^k\quad & :\,i \ge 0\\
		  d^{-i}b^jc^k\quad & :\,i<0
		  \end{array}\right.\quad
		  i \in \mathbb{Z},\>j,k \in \mathbb{N} 
\end{equation}
form a vector space basis of $\A$.

For $\lambda,\mu \in \k^\times$ there are
unique automorphisms $\sigma_{\lambda, \mu},\tau_{\lambda, \mu}$ 
of $\A$ with
$$ \sigma_{\lambda, \mu}(a)=\lambda a,\quad
	\sigma_{\lambda, \mu}(b)=\mu b,\quad
	\sigma_{\lambda, \mu}(c)=\mu^{-1} c,\quad
	\sigma_{\lambda, \mu}(d)=\lambda^{-1} d,$$
$$
	\tau_{\lambda, \mu}(a)=\lambda a,\quad
	\tau_{\lambda, \mu}(b)=\mu^{-1} c,\quad
	\tau_{\lambda, \mu}(c)=\mu b,\quad
	\tau_{\lambda, \mu}(d)=\lambda^{-1} d,
$$
and all automorphisms of $\A$ are of this form
(see \cite{joseph}). For later use, we
compute for all $ \sigma_{\lambda,\mu}$ the 
twisted commutators
\begin{eqnarray}\label{twico}
		  e_{i,j,k} a - \lambda a e_{i,j,k} 
\!\!\!\!&=&\!\!\!\! (q^{-j-k} - \lambda) e_{i+1,j,k} \nonumber\\ 
\!\!\!\!&&\!\!\!\! + \left\{
\begin{array}{ll}
0 \quad & :\,i \ge 0\\
(q^{-j-k-1}- \lambda q^{-1-2i}) e_{i+1,j+1,k+1}\quad &
:\,i<0,
\end{array}\right.\nonumber\\ 
		  e_{i,j,k} b - \mu b e_{i,j,k} 
\!\!\!\!&=&\!\!\!\! (1-\mu q^{-i}) e_{i,j+1,k},\\   
		  e_{i,j,k} c - \mu^{-1} c e_{i,j,k} 
\!\!\!\!&=&\!\!\!\! (1-\mu^{-1}q^{-i}) e_{i,j,k+1},\nonumber\\ 
		  e_{i,j,k} d - \lambda^{-1} d e_{i,j,k} 
\!\!\!\!&=&\!\!\!\!(q^{j+k} - \lambda ^{-1}) e_{i-1,j,k} \nonumber\\
\!\!\!\!&&\!\!\!\! + \left\{
\begin{array}{ll}
0 \quad & :\,i \le 0\\
(q^{j+k+1}- \lambda ^{-1} q^{1-2i}) e_{i-1,j+1,k+1}\quad &
:\,i>0 \end{array}\right..\nonumber 
\end{eqnarray} 
Finally, recall 
(see \cite{chef})
that the standard 
Hopf algebra structure on $\A$ 
admits a so-called universal r-form 
$\r : \A \otimes \A \rightarrow \k$. 
This can be used to define a braiding 
$$
		  \Psi : \A \otimes \A \rightarrow \A \otimes
		  \A,\quad
		  x \otimes y \mapsto 
		  \r(y_{(1)},x_{(1)}) y_{(2)} \otimes 
		  x_{(2)}\r(S(y_{(3)}),x_{(3)}),
$$
where $x \mapsto x_{(1)} \otimes x_{(2)}$
is the coproduct in Sweedler notation and 
$S : \A \rightarrow \A$ is the antipode.
This braiding should be
considered as a quantum analogue of the tensor flip 
and is used in the standard way to define 
\begin{eqnarray}\label{wedgep}
		  x \wedge y&:=& (\mathrm{id} - \Psi)(x \otimes y),\nonumber\\
		  x \wedge y \wedge z 
&:=& (\mathrm{id} - \Psi_{1,2}-\Psi_{2,3}+
		  \Psi_{2,3} \circ \Psi_{1,2} + \Psi_{1,2} \circ \Psi_{2,3}\\
&& -\Psi_{1,2} \circ \Psi_{2,3} \circ
		  \Psi_{1,2})(x \otimes y \otimes z)\nonumber
\end{eqnarray}
where $\Psi_{2,3}:=\mathrm{id} \otimes \Psi$
and $\Psi_{1,2}:=\Psi \otimes \mathrm{id}$.
On generators, we have
$$
\left[\begin{array}{cccc}
\r(a , a) & \r( a,b) & \r(a ,c) & \r( a,d)\cr
\r(b , a) & \r(b ,b) & \r(b ,c) & \r( b,d)\cr
\r( c, a) & \r(c ,b) & \r(c ,c) & \r(c ,d)\cr
\r(d , a) & \r(d ,b) & \r(d ,c) & \r(d ,d)\cr
 \end{array}\right]
=
q^{-1/2} \left[\begin{array}{cccc}
 q & 0 & 0 & 1\cr
  0 & 0 & 0 & 0\cr
   0 & q - q^{-1}  & 0 & 0\cr
  1 & 0 & 0 & q \cr
 \end{array}\right],
$$
where $q^{-1/2}$ is a fixed solution of $z^2=q^{-1}$, and
\begin{eqnarray}
&& \Psi (a \otimes a)=a \otimes a,\quad
		  \Psi (a \otimes b)=q^{-1}b \otimes a+
		  (1-q^{-2})a \otimes b,\nonumber\\ 
&& \Psi (a \otimes c)=qc \otimes a,\quad\Psi (a \otimes d)=
		  d \otimes a+(q-q^{-1})c \otimes b,\nonumber\\ 
&& \Psi (b \otimes a)=q^{-1}a \otimes b,\quad
		  \Psi (b \otimes b)=b \otimes b,\quad
		  \Psi (b \otimes c)=c \otimes b,\nonumber\\ 
&& \Psi (b \otimes d)=qd \otimes b,\quad
		  \Psi (c \otimes a)=
		  qa \otimes c+(1-q^2)c \otimes a,\nonumber\\ 
&& \Psi (c \otimes b)=b \otimes c
		  -(q-q^{-1})^2 c \otimes b
		  +(q-q^{-1})a \otimes d+(q^{-1}-q)d \otimes a,\nonumber\\ 
&& \Psi (c \otimes c)=c \otimes c,\quad
		  \Psi (c \otimes d)=q^{-1}d \otimes c+(1-q^{-2})c
		  \otimes d,\nonumber\\
&& \Psi (d \otimes a)=
		  a \otimes d-(q-q^{-1})c \otimes b,\quad
		  \Psi (d \otimes b)=
		  qb \otimes d+(1-q^2)d \otimes b,\nonumber\\ 
&& \Psi (d \otimes c)=q^{-1}c \otimes d,\quad
		  \Psi (d \otimes d)=d \otimes d.\nonumber 
\end{eqnarray}

 \subsection{Results}\label{hoho2} 
Here we recall from \cite{kt} the description of 
$H_n(\A,{}_\sigma \A)$ for $\A=\k_q[SL(2)]$ and
$\sigma = \sigma_{\lambda,\mu}$, but we simplify the
presentation and 
also correct some errors. Throughout, elements of
$H_n(\A,{}_\sigma \A)$ 
are represented in the normalised Hochschild complex
$\bar C_\bullet(\A,{}_\sigma \A)$.

\subsubsection{$n=0$}\label{nistnull} 
$H_0(\A,{}_\sigma \A)$ is easily computed directly, using the canonical complex 
$C_\bullet(\A,{}_\sigma \A)$. Since 
\begin{equation}\label{leubnuetz}
		  x \otimes yz=xy \otimes z+\sigma(z)x \otimes y
		  -\b(x,y,z),
\end{equation} 
the boundary operator 
$\b$ on $C_1(\A,{}_\sigma \A)$ satisfies
$$
		  \b(x,yz)=\b(xy,z)+\b(\sigma (z)x,y),
$$ 
so its image is
spanned by 
$\b(e_{i,j,k},a)$, $\b(e_{i,j,k},b)$, 
$\b(e_{i,j,k},c)$ and
$\b(e_{i,j,k},d)$,
that is, by the twisted commutators (\ref{twico}).
This yields a description of 
$H_0(\A,{}_{\sigma} \A)$ that can be
summarised in compact form as follows:
we first define 
\begin{equation}\label{brot}
\omega_{r,i}:=b^ic^{r-i},\;
\Slambda := \left\{
		  \begin{array}{ll}
		  \mathbb{N}\backslash \{
			N-2r\,:\,r \ge 1\;\} \; 
			& :\;\lambda = q^{-N}, \; N \ge 2\\
		   \mathbb{N}\; & :\; \otherwise .
		  \end{array}\right.
\end{equation}
Then the following set of 
homology classes is
a vector space basis of
$H_0(\A,{}_{\sigma} \A)$:
\begin{eqnarray}\label{hnullbasis}
&& \{[a^i],[d^i]\,:\,i \ge 0\,|\,\lambda=1\}\nonumber\\
&\cup& \{[b^j],[c^j]\,:\,j\in\Slambda\,|\,\mu=1\}\nonumber\\  
&\cup& \{[\omega_{N,i}]\,:\,1 \le i \le N-1\, | \,
		  \lambda=q^{-N},N \ge 2,\;\mu=1\}
		 \\ 
&\cup& \{[ e_{M,N,0} ],[e_{-M,0,N}]\,|\,\lambda = q^{-N}, N > 0,\;
		  \mu = q^M,M\neq 0\}.\nonumber
\end{eqnarray} 
We use the convention that $[x^0 ] := [1]$, for any $x \in \{ a,b,c,d \}$, and multiple occurrences of any $[y]$ are counted only once. 
To ensure the notation  is clear,
consider the case $\mu = 1$, 
$\lambda = q^{-N}$, $N \ge 2$ (
 Case 2 in \cite{kt}). 
Then (\ref{hnullbasis}) gives a basis 
\begin{equation}\label{mcquillan}
\{ \;[b^j ],\;[c^j ] \,:\, j \in
\Slambda \;\} \;
\cup \;\{ \;[\omega_{N,i}]\; : \; 1 \le
i \le N-1 \; \},
\end{equation} 
with $\Slambda = 
\mathbb{N}\backslash\{N-2,N-4,N-6,\ldots\}$, 
exactly as in \cite{kt}, p343.

Although  $H_0(\A,{}_\sigma \A)$ was computed  
correctly in Section 4.3 of \cite{kt}, the overview on p.328-329 therein claimed
that $H_0(\A,{}_\sigma \A)$ 
is infinite-dimensional only for $ \mu = 1$,
whereas it should correctly read for
$ \mu = 1$ or $ \lambda = 1$.

\subsubsection{$n=1$} 
We also used  
$C_\bullet(\A,{}_\sigma \A)$ to compute 
$H_1(\A,{}_\sigma \A)$.
By (\ref{leubnuetz}),
$H_1(\A,{}_\sigma \A)$ is generated by the classes of 
linear combinations of 
$e_{i,j,k} \otimes a$,
$e_{i,j,k} \otimes b$,
$e_{i,j,k} \otimes c$ and
$e_{i,j,k} \otimes d$. These are mapped by $\b$
to the twisted commutators (\ref{twico}), and it is
simple to compute which linear combinations of
these tensors defines a cycle and which of these are
homologous to each other (see \cite{kt}).  
This gives the following
vector space basis of
$H_1(\A,{}_\sigma \A)$:
\begin{eqnarray}\label{heinsbasis}
&& \{[(1-\mu^{-1}) d \otimes a + (q - q^{-1}) b \otimes c]\,|\,\lambda=1\} \nonumber\\ 
&\cup& \{[a^i \otimes a],[d^i \otimes d]\,:\,
		  i \ge 0\,|\,\lambda=1\} \nonumber\\ 
&\cup& \{[b^{j-1} \otimes b],[c^{j-1} \otimes c]\,:\,j\in\Slambda\,|\, \mu=1\}\nonumber\\ 
&\cup& \{[\omega_{N-1,i} \otimes b],
		  [\omega_{N-1,i+1} \otimes c] : 0 \le i \le N-2 | 
		  \lambda = q^{-N},N \ge 2,\mu=1\}\nonumber\\ 
&\cup& \{[a^{M-1}b^N \otimes a],
		  [a^M b^{N-1} \otimes b]\,|\,
		  \lambda=q^{-N},\mu = q^M,\;M,N>0\}\\ 
&\cup& \{[d^M c^{N-1} \otimes c],[d^{M-1}c^N \otimes d]\,|\,
		  \lambda=q^{-N},\mu = q^M,\;M,N>0\}\nonumber\\ 
&\cup& \{[d^{M-1} b^N \otimes d],
		  [d^M b^{N-1} \otimes b]\,|\,
		  \lambda=q^{-N},\mu = q^{-M},\;M,N>0\}\nonumber\\ 
&\cup& \{[a^M c^{N-1} \otimes c],[a^{M-1} c^N \otimes a]\,|\,
		  \lambda=q^{-N},\mu = q^{-M},\;M,N>0\}.\nonumber 
\end{eqnarray} 
Here $\Slambda$ and $ \omega_{r,i}$ are as in (\ref{brot}) 
and we use the formal notation (neither
$b$ nor $c$ are invertible elements)
$$
		  [c^{-1} \otimes c]:=[b \otimes c],\quad
		  [b^{-1} \otimes b]:=[c \otimes b]
$$
which appears in the above set for $\mu =1$ except when 
$ \lambda = q^{-N}$ with $N=2r$, $r>0$, but should be counted
only once: 
\begin{lemma}\label{laetst}
If $\lambda \neq q^{-2}$, then 
$[b \otimes c]=- \mu^{-1} [c \otimes b]$
 for all $\mu \in \mathbb{C}$.
\end{lemma}
\begin{pf}
 This follows from 
$[1 \otimes bc]=[b \otimes c]+\mu^{-1} [c \otimes b]$
(a special case of (\ref{leubnuetz})) together with 
$
		  \b(1 \otimes (a \otimes d-\lambda^{-1}d \otimes a+
		  (1-\lambda^{-1}) \otimes 1)) =
		  (\lambda^{-1}q^{-1}-q) \otimes bc. 
$
\end{pf}

\subsubsection{$n=2$}\label{nzwei}  
In higher degrees, working with the canonical complex is
no longer feasible, but we showed in \cite{kt}, 
Proposition~4.1 that the trival left 
$\A$-module $\k$ (on which $\A$ acts by the counit 
$ \varepsilon $ of the standard Hopf algebra structure)
admits a noncommutative 
Koszul resolution of the form
\begin{equation}\label{koszul}
\xymatrix
{0 \ar[r] & \A \ar[r] 
& \A^3 \ar[r]
& \A^3 \ar[r] 
& \A \ar[r] 
& 0}
\end{equation} 
with morphisms given by the matrices 
$$
		  \left(
		  \begin{array}{ccc}
		  c  & -b  & q^{-2}a-1 
		  \end{array}
		  \right),\quad
		  \left(
		  \begin{array}{ccc}
			b & 1-q^{-1}a & 0\\
			c & 0 & 1-q^{-1}a\\
			0 & c & -b\end{array}
			\right),\quad
			\left(
			\begin{array}{c}
			a-1 \\ b \\ c 
			\end{array}
			\right)
$$
that operate by right multiplication on row vectors.
Then we computed 
$H_n(\A,{}_\sigma \A)$ for $n \ge 2$, 
using the alternative derived functor 
description of $H_\bullet(\A,\M)$ for Hopf algebras as in
Feng and Tsygan \cite{ft}.

The result is that 
$H_2(\A,{}_\sigma \A)=0$ except when 
$\lambda=q^{-N}$, $N>0$, and in this
case one has
$$
		  \mathrm{dim}_\k H_2(\A,{}_\sigma
		  \A) =
		  \left\{\begin{array}{ll}
		  2(N-1)\quad & :\,\mu=1\\
			2 \quad & :\,\mu=q^{\pm M},M>0\\
			0 \quad & :\,\mu \notin q^{\mathbb{Z}}\\
			 \end{array}\right.
$$
A linear basis is given by
\begin{eqnarray}\label{dimH2}
&& \{[\omega_2 (N-2,i)],\;[\omega_2^{'} (N-2,i)]\, : \,
		  0 \le i \le N-2 \, | \, \mu=1\} \nonumber\\ 
&\cup& \{[a^{M-1} b^{N-1} \otimes (b \wedge a)],\;
		  [d^{M-1} c^{N-1} \otimes (d \wedge c)]\,|\,
		  \mu = q^M,M>0\}, \\
&\cup& \{[a^{M-1} c^{N-1} \otimes (a \wedge c)],\;
		  [d^{M-1} b^{N-1} \otimes (b \wedge d)]\,|\,
		  \mu = q^{-M},M>0\},\nonumber
\end{eqnarray} 
where 
$x \wedge y$ was defined in (\ref{wedgep}) and
\begin{eqnarray}
		  \omega_2 (r,i) 
&:=& \omega_{r,i}\bigl(
	bc \otimes (a \wedge d) -
	bd \otimes (a \wedge c)+ \nonumber\\ 
&&da \otimes (b \wedge c) 
	- q^{-1} ca
	 \otimes (b \wedge d)\bigr),\nonumber\\
		  \omega_2^{'} (r,i) 
&:=& \omega_{r,i} \otimes (b \wedge c).\nonumber
\end{eqnarray}
These formulae differ from
\cite{kt} by a correction we found
after S.~Launois pointed out to us an
inconsistency
between our results in \cite{kt} and
\cite{cr} (see \cite{LL}): 
the second half of the first sentence 
after Proposition~4.10 on  p.349 in
\cite{kt} is incorrect, 
the generators given 
 for  $\mu =1$ are all linearly independent in homology. 
\subsubsection{$n=3$}\label{n33} 
We have $H_3(\A,{}_\sigma \A)=0$ except when 
$\lambda=q^{-N}$, $N \ge 2$, $ \mu = 1$, and in this case
$$
		  \mathrm{dim}_\k H_3(\A,{}_{\sigma_{q^{-N},1}} \A)=N-1,
$$
with a basis given by the classes of
\begin{eqnarray}
		  \omega_3 (N-2,i) 
&:=& \omega_{N-2,i} \bigl(-qd \otimes (b \wedge a \wedge c) 
		  + c \otimes (b \wedge a \wedge
		  d)\bigr),\nonumber\\ 
&=& b^ic^{N-2-i}\bigl(d \otimes 
		  (-qb \otimes a \otimes c
		  +a \otimes b \otimes c
		  +q^2b \otimes c \otimes a \nonumber\\  
&&-a \otimes c \otimes b
		  -q^2 c \otimes b \otimes a
		  +qc \otimes a \otimes b) 
		  \nonumber\\ 
&&+ c \otimes (b \otimes a \otimes d 
		  -q^{-1} a \otimes b \otimes d 
		  -b \otimes d \otimes a \nonumber\\ 
&& -(q-q^{-1})b \otimes c \otimes b
		  + a \otimes d \otimes b 
		  +qd \otimes b \otimes a 
		  -d \otimes a \otimes b)\bigr),\nonumber 
\end{eqnarray}
for $0 \leq i \leq N-2$. We abbreviate
\begin{equation}\label{defda}
		  \dA := [\omega_3 (0,0)].
\end{equation}  
We remark that in the
explicit formula for $b \wedge a \wedge d$
given in \cite{kt} the last term
$-(q-q^{-1}) b^ic^{N-1-i} \otimes b \otimes c \otimes b$ 
was missing, and that the above defines
a cycle in the unnormalised complex 
$C_\bullet(\A,{}_{\sigma_{q^{-N},1}} \A)$, not just in 
$\bar C_\bullet(\A,{}_{\sigma_{q^{-N},1}} \A)$.

\subsubsection{$n>3$} 
Since the resolution (\ref{koszul}) has length 3,
$H_n(\A,{}_\sigma \A)=0$ for $n>3$.

\section{Hochschild cohomology}\label{como}
\label{twistedcocycles}
\subsection{Background}\label{nuda} 
Section~\ref{nuda} recalls 
products and Poincar\'e duality 
in Hoch\-schild (co)homology 
(see e.g.~\cite{CE,nielsuli,nt,vdb}),
as well as
twisted derivations of $\k_q[SL(2)]$ that arise from
the left and right actions of its Hopf 
dual.
 
\subsubsection{Hochschild cohomology} 
Let $\A$ be a unital associative $\k$-algebra. 
The Hochschild cohomology 
$H^n(\A,\M):=\mathrm{Ext}^n_{\A^e}(\A,\M)$
of $\A$ with coefficients in an $\A$-bimodule 
$\M$ can be computed as the
cohomology of the cochain complex
$C^\bullet(\A,\M)$ of $\k$-linear 
maps 
$ \varphi : \A^{\otimes \bullet} \rightarrow \M$ with
coboundary map given by
\begin{eqnarray}
		  (\b \varphi)(a_1,\ldots,a_{n+1}) 
&:=& a_1 \triangleright \varphi (a_2,\ldots,a_{n+1})\nonumber\\
&&+ \sum_{j=1}^n (-1)^j \varphi( a_1 , \ldots , a_j a_{j+1} , \ldots , a_{n+1} )
\nonumber\\ 
&&+(-1)^{n+1} \varphi (a_1,\ldots,a_n) \triangleleft a_{n+1}.
\nonumber
\end{eqnarray}   
This presentation of cocycles yields the
standard identification 
$$
		  H^0(\A,\M) \simeq \{m \in \M\,|\,a \triangleright m=
		  m \triangleleft a \mbox{ for all } a\},
$$
and of 
$H^1(\A,\M)$ with the space of 
derivations 
\begin{equation}\label{leibnizrule}
	\partial : \A \rightarrow \M,\quad
		  \partial (xy)=
		  x \triangleright \partial (y) + \partial (x)
		  \triangleleft y
\end{equation} 
modulo inner derivations (those of the form 
$x \mapsto x \triangleright m-m \triangleleft x$,
$m \in \M$). 

\subsubsection{The cup product} 
The cup product 
$$
	\smallsmile \,: H^m(\A,\M) \otimes H^n(\A,\N) \rightarrow
	H^{m+n}(\A,\M \otimes_\A \N)
$$
is defined on the level of cochains by
$$
		  (\varphi \smallsmile \psi) 
		  (a_1,\ldots,a_{m+n}):=
		  \varphi (a_1,\ldots,a_m) \otimes_\A 
		  \psi (a_{m+1},\ldots,a_{m+n}),
$$
where $\varphi \in C^m(\A,\M),\psi \in C^n(\A,\N)$. Since 
$$
		  \b(\varphi \smallsmile \psi) =
		  (\b \varphi) \smallsmile \psi +
		  (-1)^m \varphi \smallsmile
		  (\b \psi)
$$ 
the cup product is well-defined on the level
of cohomology. 
As a special case, we obtain for  
$ \sigma , \tau \in \mathrm{Aut} (\A)$ a map
\begin{equation}\label{spinne2}
		  H^m(\A,{}_\sigma \A) \otimes 
		  H^n(\A,{}_\tau \A) \rightarrow 
		  H^{m+n}(\A,{}_{\tau \circ \sigma} \A)
\end{equation} 
given on cochains by
\begin{equation}\label{verz}
		  (\varphi \smallsmile \psi) 
		  (a_1,\ldots,a_{m+n})=
		  \tau (\varphi (a_1,\ldots,a_m)) 
		  \psi (a_{m+1},\ldots,a_{m+n}).
\end{equation} 
Thus for any monoid 
$G \subset \mathrm{Aut} (\A)$
we obtain an $\mathbb{N} \times G$-graded 
algebra
$$
		  \Lambda^\bullet_G(\A):=
		  \bigoplus_{n \in \mathbb{N},\sigma \in G}
		  H^n(\A,{}_\sigma \A).
$$ 
We call 
the subalgebra $ \Lambda^0_G(\A)$ the $G$-twisted
centre of $\A$ and the elements of the 
$ \Lambda^0_G(\A)$-bimodule 
$ \Lambda^1_G(\A)$ (or rather the
cocycles representing them) the $G$-twisted derivations of $\A$.
Obviously, $G$ could be replaced by any monoid of
bimodules, but we will not need this in the present paper.

Note that in degree 0, (\ref{verz}) reduces to
the
\emph{opposite} product of $\A$,
\begin{equation}\label{oppo}
		  x \smallsmile y=\tau(x)y = yx,\quad
		  x \in H^0(\A,{}_\sigma \A),\quad
		  y \in H^0(\A,{}_\tau \A),
\end{equation} 
and that for $z \in H^0(\A,{}_\sigma \A)$ and 
$\partial \in C^1(\A,{}_\tau \A)$ we have
\begin{eqnarray}\label{drospa}
&&(z \smallsmile \partial)(x)= \sigma (\partial (x))
		  z=z \partial (x),\nonumber\\ 
&& (\partial \smallsmile z)(x)= \tau(z) \partial (x),\\
& \Rightarrow & \partial \smallsmile z= \tau(z) \smallsmile
			\partial.\nonumber
\end{eqnarray} 
Finally, for twisted derivations
$\partial \in C^1(\A,{}_\sigma \A),
\partial' \in C^1(\A,{}_\tau \A)$ we have
\begin{eqnarray}
&& (\partial \smallsmile \partial' +
		  (\sigma ^{-1} \circ \partial ' \circ \sigma)
		  \smallsmile \partial)(x,y) \nonumber\\ 
&=& \tau (\partial (x)) \partial ' (y) 
		  +\partial ' (\sigma (x)) \partial (y) \nonumber\\ 
&=& \partial ' (\partial (x)y) -
		  \partial ' (\partial (x))y+
		  \partial ' (\sigma (x) \partial (y))-
		  \tau (\sigma (x)) \partial '
		  (\partial (y))
		  \nonumber\\ 
&=& \partial '(\partial (xy)) -
		  \partial ' (\partial (x))y-
		  \tau (\sigma (x)) \partial ' (\partial (y))
		  \nonumber\\ 
&=& -\b(\partial ' \circ \partial)(x,y),\nonumber
\end{eqnarray}
so in cohomology 
\begin{equation}\label{eigentor}
		  [\partial] \smallsmile [\partial']=
		  -[\sigma ^{-1} \circ \partial ' \circ \sigma]
		  \smallsmile [\partial] \in 
		  H^2(\A,{}_{\tau \circ \sigma} \A).
\end{equation} 

\subsubsection{The cap product}\label{duality33} 
The duality between Hochschild homology and cohomology
results from the cap product pairing
$$
		  \smallfrown \, : H_n(\A,\M) \otimes H^m(\A,\N) \rightarrow
		  H_{n-m} (\A,\M \otimes_\A \N),\quad m \le n
$$
defined on (co)chains  by evaluation,
$$
		   (a_0 \otimes \ldots \otimes
		  a_n) \smallfrown \varphi=
		  a_0 \otimes_\A \varphi(a_1,\ldots,a_m) \otimes 
		  a_{m+1} \otimes \ldots \otimes a_n. 
$$
For 
$m=n$, the pairing $\smallfrown$ becomes the duality
pairing from 
\cite{loday}, Section~{1.5.9} after identifying 
$$
H_0(\A,\M \otimes_\A \N )=
\M \otimes_\A \N \otimes_{\A^e} \A \simeq
\M \otimes_{\A^e} \N.
$$
Taking 
$\N=\M^*=\mathrm{Hom}_\k(\M,\k)$ and composing with the
canonical 
evaluation map $\M \otimes_{\A^e} \M^* \rightarrow \k$ 
gives the duality
pairing 
$$
		  H_n(\A,\M) \otimes H^n(\A,\M^*) \rightarrow
		  \k.
$$
By the universal coefficient theorem, this 
yields an isomorphism 
$H^n(\A,\M^*) \simeq (H_n(\A,\M))^*$. In this way
a Hochschild cocycle 
$ \varphi \in C^n(\A,\M^*)$ will
usually be viewed as a $\k$-linear map
$\M \otimes \A^{\otimes n} \rightarrow \k$.

For any $G \subset \mathrm{Aut} (\A)$, the cap product 
endows 
$$
		  \Omega_\bullet^G(\A):=
		  \bigoplus_{n \in \mathbb{N},\sigma \in G}
		  H_n(\A,{}_{\sigma^{-1}} \A)
$$
with the structure of 
an $ \mathbb{N} \times G$-graded 
(right) module over $ \Lambda^\bullet_G(\A)$
(here we use the convention that a
homogeneous element of a graded ring lowers
the degree of an element of a homogeneous module by 
its degree).
Remember to take into account the
identification
${}_\sigma \A \otimes_\A {}_\tau \A \rightarrow 
{}_{\tau \circ \sigma} \A$ :
explicitly, the action of 
$ \varphi \in C^m(\A,{}_\tau \A)$ on 
$a_0 \otimes \ldots \otimes a_n \in C_n(\A,{}_\sigma \A)$
is given by
$$
		  (a_0 \otimes \ldots \otimes a_n) \smallfrown
		  \varphi =
		  \tau (a_0) \varphi (a_1,\ldots,a_m) \otimes 
		  a_{m+1} \otimes \ldots \otimes a_n \in
		  C_{n-m}(\A,{}_{\tau \circ \sigma} \A). 
$$
In particular, the cap product
with a twisted central
element 
$z \in H^0(\A,{}_\sigma \A)$ is simply given by
multiplication from the left,
\begin{equation}\label{spaetzle}
		  (a_0 \otimes \ldots \otimes a_n) \smallfrown z=
		  \sigma (a_0)z \otimes \ldots \otimes a_n=
		  za_0 \otimes \ldots \otimes a_n.
\end{equation}

The cap product is also the source of Poincar\'e-type
dualities in Hochschild (co)homology. In
particular, we have for
$\A=\k_q[SL(2)]$:
\begin{thm}\cite{cr,nielsuli}\label{pds}
For any bimodule $\N$ over 
$\A:=\k_q[SL(2)]$, the map 
\begin{equation}\label{poincare}
		  \cdot \smallfrown \dA : 
			H^n(\A,\N) \rightarrow H_{3-n}(\A,{}_{\sigma_{q^{-2},1}}\N)
\end{equation} 
is an isomorphism of $\k$-vector spaces.
\end{thm}

\subsubsection{Twisted primitive elements in the 
Hopf dual of $\k_q[SL(2)]$}\label{tutor}
For the rest of Section~\ref{como}, 
we focus on $\A=\k_q[SL(2)]$ and 
recall (see e.g.~\cite{istvan,chef} and
references therein) that the Hopf dual of the standard
Hopf algebra structure on $\A$ contains a Hopf subalgebra
$\U$ with generators
$H,K,K^{-1},E,F$ and relations 
\begin{eqnarray}
&& KK^{-1}=K^{-1}K=1,\quad
	KEK^{-1}=q^2E,\quad
	KFK^{-1}=q^{-2}F,\nonumber\\ 
&& [E,F]=\frac{K-K^{-1}}{q-q^{-1}},\quad [H,K]=0,\quad
		  [H,E]=2E,\quad
		  [H,F]=-2F.\nonumber
\end{eqnarray} 
This is the standard Drinfeld-Jimbo quantised
universal enveloping algebra 
$U_q(\mathfrak{sl}(2))$ extended by the unquantised
functional $H$ (so 
when working with formal deformations we would have 
$K=e^{\hbar H}$, $q=e^\hbar$).

The dual pairing of $\U$ and $\A$ gives two
commuting left and right actions of $\U$ on $\A$, and 
the operators assigned by these actions 
to $H,EK^{-1}$ and $F$ are (twisted) derivations
which we denote by
\begin{eqnarray}
&& {\partial^+_H}  : a,b,c,d \mapsto 
		  -a,b,-c,d,\quad
		  {\partial^-_H} : a,b,c,d \mapsto 
		  -a,-b,c,d \in C^1(\A,\A),\nonumber\\
&& {\partial^+_E} : a,b,c,d \mapsto
	qb,0,qd,0,\quad {\partial^+_F}  : a,b,c,d \mapsto
	0,a,0,c \in C^1(\A,{}_{\sigma_{q,q^{-1}}}\A),\nonumber\\ 
&& {\partial^-_E}: a,b,c,d \mapsto
	0,0,q^{-1}a,q^{-1}b,\quad {\partial^-_F} : a,b,c,d \mapsto
	c,d,0,0 \in 
	C^1(\A,{}_{\sigma_{q,q}}\A).\nonumber 
\end{eqnarray}
Due to the Leibniz rule (\ref{leibnizrule}),
this determines the derivations uniquely. 

\subsection{Results}
Here we compute the twisted
centre and twisted derivations of $\k_q[SL(2)]$. Then
we show how their cap product action on the twisted 
Hochschild homology groups can be used to determine the
homology class of a given 2- or 3-cycle.
 
\subsubsection{Twisted central elements} 
The twisted centre 
$\Lambda^0:=\Lambda^0_{\mathrm{Aut}(\A)}(\A)$ 
of $\A$ is a 
(commutative) polynomial ring in two
indeterminates: 
\begin{lemma}\label{losta}
There is an isomorphism
of graded algebras  
$$
		  \Lambda^0=\bigoplus_{N \ge 0} 
		  H^0(\A,{}_{\sigma_{q^{-N},1}} \A)
		  \simeq \k[b,c].
$$
\end{lemma}
\begin{pf}
By Poincar\'e duality (Theorem~\ref{pds}) and the
 computation of $H_3(\A,{}_\sigma \A)$ 
recalled in Section~\ref{n33}, $H^0(\A,{}_\sigma \A)$ 
vanishes except when 
$ \sigma = \sigma_{q^{-N},1}$, 
$N \ge 0$, and in this case 
it has dimension $N+1$ over $\k$. 
The monomials $\omega_{N,i}=b^ic^{N-i}$,
$0 \le i  \le N$, are $N+1$ linearly independent 
elements of $H^0(\A,{}_{\sigma_{q^{-N},1}} \A)$
and hence form a vector space basis. 
Finally, (\ref{oppo}) gives 
$\omega_{r,i} \smallsmile \omega_{s,j}=
\omega_{r,i} \omega_{s,j} = \omega_{r+s,i+j}$ since $b$ and $c$ commute.  
\end{pf}

The cap product action of $ \Lambda^0$
gives additional structure 
to the results of our computations 
of $H_\bullet(\A,{}_\sigma \A)$. For example, the
direct sum of the nontrivial
$H_3(\A,{}_\sigma \A)$ forms 
a free module of rank
one over $\Lambda^0$ with generator $\dA$.
Similarly we will see below that applying twisted
derivations to $\dA$ leads for example to the
2-dimensional 
$H_2(\A,{}_{\sigma_{q^{-N},q^{\pm M}}}\A)$ in (\ref{dimH2}).
 
\subsubsection{Detecting nontrivial 2-cycles}\label{fest}
For large $n$ and $m$ it is
typically difficult to decide whether a given cycle
$c \in C_n(\A,{}_\sigma \A)$ and cocycle
$\varphi \in C^m(\A,{}_\tau \A)$ have 
nontrivial classes in (co)homology. A sufficient
criterion is that  
$c \smallfrown \varphi \in C_{n-m}(\A,{}_{\tau \circ \sigma} \A)$
is nontrivial in homology, and this may be easy to
verify for small $n-m$.
We now give an example of this kind 
for $\A=\k_q[SL(2)]$, whose result 
is used later in
the computation of twisted cyclic homology and also 
to determine $\Lambda^1_{\mathrm{Aut} (\A)}(\A)$. 
\begin{lemma}\label{tester1}
Abbreviate 
$\partial:=
		  \frac{1}{2} 
		  ({\partial^+_H}+{\partial^-_H})$ and
$\partial':=
		  -{\partial^-_H}$.
Then: 
 $$
		  [\omega_2(N-2,i)] \smallfrown [\partial]
		  =[\omega'_2(N-2,i)] \smallfrown [\partial']=
		  [\omega_{N-1,i+1} \otimes c] +
		[\omega_{N-1,i} \otimes b],
$$
$$
	  [\omega'_2(N-2,i)] \smallfrown
		  [\partial]=[\omega_2(N-2,i)] \smallfrown [\partial']=0.
$$ 
\end{lemma}
 \begin{pf}
For
$\sigma = \sigma_{q^{-2},1}$ we have
in $\bar C_1(\A,{}_\sigma \A)$:
\begin{eqnarray}
&& \b(bc \otimes a \otimes d)=
		  q^{-2}abc \otimes d-qbc \otimes bc
		  +q^2dbc \otimes a,\nonumber\\  
&& \b(bc \otimes b \otimes c)=
		  b^2c \otimes c-bc \otimes bc+bc^2 \otimes b,
		  \nonumber\\ 
&& \b(ca \otimes (d \wedge b))=(q^3-q)bc^2
		  \otimes b,\nonumber\\ 
&& \b(ba \otimes (d \wedge c))=(q-q^{-1})b^2c
		  \otimes c.\nonumber
\end{eqnarray} 
Using this, we compute directly that
\begin{eqnarray}
		  [\omega_2(0,0)] \smallfrown [{\partial^+_H}] 
&=& 2[-q^{-2}abc \otimes d-
		  q^{2}dbc \otimes a+qbc^2 \otimes b
		  \nonumber\\ 
&& +q^{-1}b^2c \otimes c+
		  c \otimes b+b \otimes c]\nonumber\\ 
&=&  2 (q^{-1}-q) [ \omega_{3,2} \otimes c]+
		   2 [\omega_{1,1} \otimes c] +
			2 [\omega_{1,0} \otimes b]\nonumber\\ 
&=& 2 [\omega_{1,1} \otimes c] +
			2 [\omega_{1,0} \otimes b],\nonumber
\end{eqnarray}
that $[\omega_2(0,0)] \smallfrown [{\partial^-_H}] 
= 0$, and that  
$$
		  [\omega'_2(0,0)] \smallfrown [{\partial^+_H}]
		  =-[\omega'_2(0,0)] \smallfrown [{\partial^-_H}]  
		  = [c \otimes b+b \otimes c].
$$
The claim follows by
$\Lambda^0_{\mathrm{Aut} (\A)}(\A)$-linearity of the products.
\end{pf}
\begin{cor}\label{swird}
The $2(N+1)$ cohomology classes 
$$
		  [\omega_{N,i}] \smallsmile [\partial],\quad
		  [\omega_{N,i}] \smallsmile [\partial'],\quad
		  0 \le i \le N,
$$ 
are linearly independent in 
$H^1(\A,{}_{\sigma_{q^{-N},1}}\A)$.
\end{cor}
 
\subsubsection{Twisted derivations}\label{huhu1} 
We now describe the twisted derivations (modulo inner
derivations)
$ \Lambda^1:=\Lambda^1_{\mathrm{Aut}(\A)}(\A)$ 
as a bimodule over the twisted centre 
$ \Lambda^0 $, and study their relations under
the cup product.

First, note that for all $i>0$ the cochains
\begin{eqnarray}
&& C^1(\A,{}_{\sigma_{q,q^{-i}}}\A) \ni
	{\partial^+_E}\smallsmile d^{i-1} : a,b,c,d \mapsto
	qd^{i-1}b,0,qd^i,0,\nonumber\\ 
&& C^1(\A,{}_{\sigma_{q,q^{-i}}}\A) \ni
	{\partial^+_F} \smallsmile a^{i-1} : a,b,c,d \mapsto
	0,a^i,0,a^{i-1}c,\nonumber\\ 
 && C^1(\A,{}_{\sigma_{q,q^i}}\A) \ni
	{\partial^-_E}\smallsmile a^{i-1} : a,b,c,d \mapsto
	0,0,q^{-1}a^i,q^{-1}a^{i-1}b,\nonumber\\ 
&& C^1(\A,{}_{\sigma_{q,q^i}}\A) \ni
	{\partial^-_F} \smallsmile d^{i-1} : a,b,c,d \mapsto
	d^{i-1}c,d^i,0,0,\nonumber 
\end{eqnarray}
are (twisted) derivations, although 
$a^{i-1},d^{i-1} \notin \Lambda^0$. The point is that
for example $\partial^+_E(\A) \subset \B$, where 
$\B \subset \A$ is the subalgebra generated 
by $b,d$, and we have $d^{i-1} \in 
\Lambda^0_{\mathrm{Aut} (\B)}(\B)$ with twisting
automorphism extending to the whole of $\A$.
This implies the claim, which of course can also be
verified directly.
We denote the classes of these derivations 
in cohomology by 
$\partial_i^+,\partial_{-i}^+,\partial_i^-,\partial_{-i}^-$,
respectively, and define
$\partial_0^\pm:=[\partial^\pm_H]$.  

\begin{lemma}\label{derifs}
As a left $\Lambda^0$-module, 
$ \Lambda^1$ is generated by  
$\{\partial_i^\pm,i \in \mathbb{Z}\}$, and
\begin{eqnarray}\label{dro}
&& [b] \smallsmile \partial^\pm_{-i}=0,\quad 
  		  [c] \smallsmile \partial^\pm_i=0,\quad i >
  		  1,\nonumber\\  
&&
		  \partial_i^\pm \smallsmile [b^jc^k] =
			q^{\pm i(j-k)} [b^jc^k] \smallsmile
			\partial_i^\pm,\quad i \in \mathbb{Z},\>j,k \in
			\mathbb{N},\\
&&		  \partial_i^{\varepsilon} \smallsmile \partial_j^\delta =
		  -q^{\delta j + \varepsilon |i| \mathrm{sgn}(j)} 
		  \partial_j^\delta
		  \smallsmile \partial_i^{\varepsilon},\quad
		  i,j \in \mathbb{Z},\>\varepsilon,\delta \in \{-1,+1\}.
		  \nonumber
\end{eqnarray} 
\end{lemma} 
\begin{pf}
The second and third relations in 
(\ref{dro}) are
computed directly from (\ref{drospa}) and
(\ref{eigentor}). Next, 
Section~\ref{nzwei} and
(\ref{poincare}) give 
$$
		  \mathrm{dim}_\k 
		  H^1(\A,{}_{\sigma_{\lambda,\mu}}\A) = 
		  \left\{\begin{array}{ll}
		  2(N+1)\quad & :\,\lambda=q^{-N},N>0,\;\mu=1,\\
		  2 \quad & :\,\lambda=q^{-N},N>0,\;\mu=q^{M},M \neq 0,\\
		  0 \quad & :\,\otherwise. \end{array}\right.
$$ 
Corollary~\ref{swird} implies that 
$\{[\omega_{N,i}] \smallsmile \partial^\pm_0\,|\,0 \le i \le N\}$ 
is a linearly independent subset of 
$H^1(\A,{}_{\sigma_{q^{-N},1}}\A)$, so it is a basis
for dimension reasons. 
The $\partial^\pm_i$, $i \neq 0$, are treated similarly. 
For example, we can use that
for $i,j>0$ 
$$
		  [\omega_2'(0,0)] \smallfrown 
		  (\partial^\pm_i \smallsmile
		  \partial^\pm_{-j})=-q^{\pm(1-i)}[e_{\pm(j-i),0,0}]
		  \in H_0(\A,{}_{\sigma_{1,q^{\mp(i+j)}}} \A),
$$
and that by the third line in (\ref{drospa}),
$$
		  \partial_i^{\varepsilon} \smallsmile \partial_j^\delta =
		  q^{(\mathrm{sgn}(i)+\mathrm{sgn}(j))
		  (\varepsilon |i|+\delta |j|)}
		  \partial_i^{\varepsilon} \smallsmile \partial_j^\delta,
$$
so $\partial_i^{\varepsilon} \smallsmile
\partial_j^\delta =0$
except when
$\mathrm{sgn}(j)=-\mathrm{sgn}(i)$ or 
$\delta=-\varepsilon,|j|=|i|$. 
\end{pf}

\subsubsection{Detecting nontrivial 3-cycles}\label{twitra}
We present computations similar 
to Section~\ref{fest},  this time acting on 
$\omega_3(0,0)$ to decide whether 
3-cycles are nontrivial. 
\begin{lemma}\label{tester2}
In $H_1(\A,{}_{\sigma_{q^{-N},1}} \A)$, we have
$$
		  [\omega_3(N-2,i)] \smallfrown
		  ([{\partial^+_H}] \smallsmile [{\partial^-_H}])=
		  2([\omega_{N-1,i+1} \otimes c] +
		[\omega_{N-1,i} \otimes b]).
$$
In particular, 
$\cdot \smallfrown
([{\partial^+_H}] \smallsmile [{\partial^-_H}]) :
\Omega_3^{\mathrm{Aut} (\A)}(\A)
\rightarrow \Omega_1^{\mathrm{Aut} (\A)}(\A)$
is injective.
\end{lemma}
\begin{pf}
Directly, 
$
		  [\omega_3(0,0)] \smallfrown 
		  [\partial^-_H]=
		  -[\omega_2(0,0)]
$ and 
$\omega_3(r,i)=\omega_3(0,0) \smallfrown
 \omega_{r,i}$. The result then follows from Lemmas~\ref{tester1}
and ~\ref{derifs}.
\end{pf}

As we shall see below, this Lemma is  strong
enough to detect nontrivial homology classes. However,
it would be  easier (and more standard) to apply a
further twisted derivation to obtain 0-cycles whose
classes in homology are even simpler to control than
those of 1-cycles. While this works  for
the fundamental class $\dA=[\omega_3(0,0)]$, 
it is unfortunately not possible for all the
other generators $[\omega_3(r,i)]$, $r>0$ : the orbit of 
$\dA$ under the cap product action of $ \Lambda^1$
is determined completely by Lemma~\ref{derifs} and 
the 
 relations
\begin{eqnarray}
&& \omega_3(0,0) \smallfrown ({\partial^+_H} \smallsmile 
		  {\partial^+_E} \smallsmile 
		  {\partial^+_F}) = 
		  (q^{-1}-q)bc+1,
		  \nonumber\\ 
&& \omega_3(0,0) \smallfrown ({\partial^+_H} \smallsmile 
		  {\partial^+_E} \smallsmile 
		  {\partial^-_H}) = 
		  2(q^4-1)db^2c-2qdb,
		  \nonumber\\ 
&& \omega_3(0,0) \smallfrown ({\partial^+_H} \smallsmile 
		  {\partial^+_E} \smallsmile 
		  {\partial^-_E}) = 
		  (q-q^{-3})b^3c-q^{-2}b^2,
		  \nonumber\\ 
&& \omega_3(0,0) \smallfrown ({\partial^+_H} \smallsmile 
		  {\partial^+_E} \smallsmile 
		  {\partial^-_F}) = (q-q^5)d^2bc+d^2,
		  \nonumber\\ 
&& \omega_3(0,0) \smallfrown ({\partial^+_H} \smallsmile 
		  {\partial^+_F} \smallsmile 
		  {\partial^-_H}) = 
		  2(q-q^{-3})abc^2-2ac,
		  \nonumber\\ 
&& \omega_3(0,0) \smallfrown ({\partial^+_H} \smallsmile 
		  {\partial^+_F} \smallsmile 
		  {\partial^-_E}) = 
		  (q^{-1}-q^{-5})a^2bc-a^2,
		  \nonumber\\ 
&& \omega_3(0,0) \smallfrown ({\partial^+_H} \smallsmile 
		  {\partial^+_F} \smallsmile 
		  {\partial^-_F}) = 
		  (q^{-1}-q^3)bc^3+(2-q^2)c^2,
		  \nonumber\\ 
&& \omega_3(0,0) \smallfrown ({\partial^+_H} \smallsmile 
		  {\partial^-_H} \smallsmile 
		  {\partial^-_E}) = 
		  2(q^{-4}-q^{-2})ab^2c+2q^{-1}ab,
		  \nonumber\\ 
&& \omega_3(0,0) \smallfrown ({\partial^+_H} \smallsmile 
		  {\partial^-_H} \smallsmile 
		  {\partial^-_F}) = 2(q-q^3)dbc^2+2dc,
		  \nonumber\\ 
&& \omega_3(0,0) \smallfrown ({\partial^+_H} \smallsmile 
		  {\partial^-_E} \smallsmile 
		  {\partial^-_F}) = 2(1-q^2)b^2c^2+2q^{-1}bc+1,
		  \nonumber\\ 
&& \omega_3(0,0) \smallfrown ({\partial^+_E} \smallsmile 
		  {\partial^+_F} \smallsmile 
		  {\partial^-_H}) = 
		  2(q^{-4}-q^2)b^2c^2+(2q^{-3}-q+q^{-1})bc+1,
		  \nonumber\\ 
&& \omega_3(0,0) \smallfrown ({\partial^+_E} \smallsmile 
		  {\partial^+_F} \smallsmile 
		  {\partial^-_E}) = 
		  (q^{-7}-q^{-1})ab^2c+q^{-4}ab,
		  \nonumber\\ 
&& \omega_3(0,0) \smallfrown ({\partial^+_E} \smallsmile 
		  {\partial^+_F} \smallsmile 
		  {\partial^-_F}) = 
		  (q^4-q^{-2})dbc^2-q^{-1}dc,
		  \nonumber\\ 
&& \omega_3(0,0) \smallfrown ({\partial^+_E} \smallsmile 
		  {\partial^-_H} \smallsmile 
		  {\partial^-_E}) = 
		  (q-q^{-3})b^3c-q^{-2}b^2,
		  \nonumber\\ 
&& \omega_3(0,0) \smallfrown ({\partial^+_E} \smallsmile 
		  {\partial^-_H} \smallsmile 
		  {\partial^-_F}) = 
		  (q^5-q)d^2bc-d^2,
		  \nonumber\\ 
&& \omega_3(0,0) \smallfrown ({\partial^+_E} \smallsmile 
		  {\partial^-_E} \smallsmile 
		  {\partial^-_F}) = 
		  (q^5-q)db^2c-db,
		  \nonumber\\ 
&& \omega_3(0,0) \smallfrown ({\partial^+_F} \smallsmile 
		  {\partial^-_H} \smallsmile 
		  {\partial^-_E}) = 
		  (q^{-5}-q^{-1})a^2bc+a^2,
		  \nonumber\\ 
&& \omega_3(0,0) \smallfrown ({\partial^+_F} \smallsmile 
		  {\partial^-_H} \smallsmile 
		  {\partial^-_F}) = 
		  (q^{-1}-q^3)bc^3+(2-q^2)c^2,
		  \nonumber\\ 
&& \omega_3(0,0) \smallfrown ({\partial^+_F} \smallsmile 
		  {\partial^-_E} \smallsmile 
		  {\partial^-_F}) = 
		  (q^{-2}-q^2)abc^2+q^{-1}ac,
		  \nonumber\\ 
&& \omega_3(0,0) \smallfrown ({\partial^-_H} \smallsmile 
		  {\partial^-_E} \smallsmile 
		  {\partial^-_F}) = 1,\nonumber
\end{eqnarray}
which hold on the level of chains in the normalised
Hochschild complex. By inspection we now obtain:
\begin{lemma}\label{bloedel}
We have 
\begin{eqnarray}
&& \dA \smallfrown 
		  [{\partial^+_H} \smallsmile  
		  {\partial^+_E} \smallsmile
		  {\partial^+_F}]=[1] \in
		  H_0(\A,{}_{\sigma_{1,q^{-2}}} \A),\nonumber\\
&& [\omega_3(r,r)] \smallfrown 
		  [{\partial^+_H} \smallsmile  
		  {\partial^-_E} \smallsmile
		  {\partial^+_E}]=[b^{r+2}] \in
		  H_0(\A,{}_{\sigma_{q^{-r},1}} \A),\nonumber\\ 
&& [\omega_3(r,0)] \smallfrown 
		  [{\partial^-_H} \smallsmile  
		  {\partial^-_F} \smallsmile
		  {\partial^+_F}]=[c^{r+2}] \in
		  H_0(\A,{}_{\sigma_{q^{-r},1}} \A),\nonumber
\end{eqnarray} 
but for all 
$  \partial_1,\partial_2,\partial_3 \in \Lambda^1$, 
$z \in \Lambda^0$, and $0<i<r$ we have
$$
		  [\omega_3(r,i)] \smallfrown 
		  [z \smallsmile \partial_1 \smallsmile 
		  \partial_2 \smallsmile \partial_3]= 0.
$$
\end{lemma}
 \begin{pf}
The first three statements are obtained
by direct computation using the above
description of the orbit of $\omega_3(0,0)$
under the cap product action of the twisted
derivations $\{{\partial^\pm_H},
{\partial^\pm_E},
{\partial^\pm_F}\}$, the fact that 
$\omega_3(r,i)=\omega_3(0,0) \smallfrown
  \omega_{r,i}$ and the commutation
  relations (\ref{dro}). 
  For the final statement, by
Lemma~\ref{derifs} we can assume 
without loss of generality 
that $\partial_1,\partial_2,\partial_3 
\in \{{\partial^\pm_H},
{\partial^\pm_E},
{\partial^\pm_F}\}$. The same Lemma
implies that we then have
\begin{eqnarray}
&& [\omega_3(r,i)] \smallfrown [z
 \smallsmile \partial_1 \smallsmile 
		  \partial_2 \smallsmile
		  \partial_3] \nonumber\\ 
&=& [\omega_3(0,0)] \smallfrown
 [\omega_{r,i} \smallsmile z
 \smallsmile \partial_1 \smallsmile 
		  \partial_2 \smallsmile
		  \partial_3] \nonumber\\ 
&=& [\omega_3(0,0)] \smallfrown
 [\partial_1 \smallsmile 
		  \partial_2 \smallsmile
		  \partial_3 \smallsmile z'] \nonumber 
\end{eqnarray} 
for some $z' \in \Lambda^0$. 
The homology class of
$\omega_3(0,0) \smallfrown 
(\partial_1 \smallsmile 
\partial_2 \smallsmile \partial_3)$
can be read off for these 
$ \partial_i$ from the list above (using
  also the commutation relations (\ref{dro}) of 
the $ \partial_i$), and
whenever the result is nonzero,
then Lemma~\ref{derifs} gives
$[\omega_{r,i} \smallsmile \partial_1 \smallsmile 
\partial_2 \smallsmile \partial_3] = 0$ for 
$0<i<r$, which implies the claim.
\end{pf}

We can now begin the proof of Theorem \ref{themainresult2}. 
We wish to compose the action of 
$[{\partial^+_H} \smallsmile  
{\partial^+_E} \smallsmile
{\partial^+_F}]$ above with a suitable twisted trace
$$
		  \int : \A \rightarrow \k,\quad
		  \int xy=\int \sigma_{1,q^{-2}} (y)x.
$$
to obtain a numerical invariant of $\dA$. 
For $\A=\k_q[SL(2)]$, the complete list of twisted traces
can be given as follows: 
for any element $[e_{i,j,k}]$ of our basis 
(\ref{hnullbasis}) of $H_0(\A,{}_\sigma \A)$
define a linear functional $\int_{[e_{i,j,k}]}$ by
$$
		  \int_{[e_{i,j,k}]} e_{r,s,t}:=\delta_{i,r} \cdot \left\{
		  \begin{array}{ll}
			\sum_{n=0}^{\infty} 
			\delta_{s,j+n} \delta_{t,k+n} (-q)^n 
			\frac{1-q^{j+k}\lambda}{1-q^{j+k+2n}\lambda} 
			\quad & :\;i=0,\>jk=0\\
			 \delta_{s,j} \delta_{t,k}  
			 \quad & :\;\otherwise.
		  \end{array}\right.  
$$ 
These  then descend to linearly independent 
functionals on 
$H_0(\A,{}_\sigma \A)$ that are
dual to the basis (\ref{hnullbasis}),
$\int_{[e_{i,j,k}]} [e_{r,s,t}]=\delta_{i,r}
\delta_{j,s} \delta_{k,t}$ for $[e_{i,j,k}],[e_{r,s,t}] \in (\ref{hnullbasis})$.

If $ \varphi \in C^n(\A,{}_\sigma \A)$ is an
$n$-cocycle and $\int \in C^0(\A,({}_\tau \A)^\ast)$ is
a twisted trace, using
${}_\sigma \A \otimes ({}_\tau \A)^* \simeq
{}_\sigma ({}_\tau \A)^* \simeq
({}_\tau \A_\sigma)^* \simeq
({}_{\sigma^{-1} \circ \tau} \A)^*$,
$\varphi \smallsmile \int$
can be identified  
with a functional on 
$H_n(\A,{}_{\sigma^{-1} \circ \tau}\A)$.
In particular 
 Lemma~\ref{bloedel} gives:
\begin{cor}
Let $\int_{[1]} : \A \rightarrow \k$ be the 
$ \sigma_{1,q^{-2}}$-twisted trace given by
$$
		  \int_{[1]} e_{r,s,t}:=\delta_{0,r}\; 
		\delta_{0,s}\; \delta_{0,t}. 
$$ 
Then
$$
		  \int_{[1]}
		  \dA \smallfrown 
		  [{\partial^+_H} \smallsmile 
		  {\partial^+_E} \smallsmile
		  {\partial^+_F}]=1.
$$
\end{cor}

Hence the linear functional
\begin{equation}\label{losfundamentos}
		  \varphi(a_0,a_1,a_2,a_3):=
			\int_{[1]} 
			\sigma_{q^2,q^{-2}}(a_0 \; \partial^+_H(a_1))\;
			\sigma_{q,q^{-1}}(\partial^+_E(a_2))\;
			\partial^+_F(a_3)
\end{equation}
is a Hochschild 3-cocycle with a nontrivial class 
in $H^3(\A,({}_{\sigma_{q^{-2},1}} \A)^*)$ that is dual
to the fundamental class 
$\dA$ in the sense that $ \varphi(\dA)=1$.

We will complete the proof of 
Theorem \ref{themainresult2} in Section \ref{cyccohom}.

\section{Cyclic homology}\label{schnoel} 
\subsection{Background}\label{precyc} 
In Section~\ref{precyc} we recall
the definition of the twisted cyclic
homology of an algebra. 
For more background see
for example \cite{loday,weibel}. 
 
\subsubsection{Paracyclic objects and
their homology} 
Paracyclic objects \cite{getzler}  
(say in an abelian category) 
slightly generalise Connes' cyclic
objects \cite{blabla}:
\begin{dfn}
A paracyclic object is a simplical 
object $(C_\bullet,\b_\bullet,\s_\bullet)$ equipped
with morphisms 
$\t : C_n \rightarrow C_n$ that satisfy (on 
$C_n$)
$$
	\b_i \t=-\t \b_{i-1},\> \s_i \t=-\t \s_{i-1},\>
	\b_0 \t=(-1)^n \b_n,\> \s_0 \t=(-1)^n \t^2 \s_n,
		  \> 1 \le i \le n. 
$$
\end{dfn}
The difference with cyclic
objects is that 
$\T:=\t^{n+1}$
is not required to be the identity $\mathrm{id}$.
However, it can be directly verified that
$\T$ commutes with all the paracyclic generators
$\t,\b_i,\s_j$. As a consequence, a cyclic object can be attached to  
any paracyclic object by passing to the coinvariants 
$C/ \mathrm{im}\, (\mathrm{id} - \T)$ of $\T$. In well-behaved
cases, there is no loss of homological information 
in this step - for example, we have
(\cite{kt}, Proposition~{2.1}):
\begin{lemma}\label{zahn}
If
$C$ is a paracyclic 
$\k$-vector space and $\T$ is diagonalisable, then
$(C,\b) \rightarrow (C/ \mathrm{im}\, (\mathrm{id}-\T),\b)$ is a
quasi-isomorphism.
\end{lemma}

Just as for cyclic objects,  
for any paracyclic object define
\begin{equation}
		  \NN:=\sum_{i=0}^{n} \t^i,\quad
			\s := (-1)^{n+1} \t \s_n,\quad
			\BB := (\mathrm{id}-\t)\s\NN,\nonumber
\end{equation} 
all acting on $C_n$, and as in the cyclic case
we have
\begin{equation}\label{skeip}
		  \b(\mathrm{id}-\t)=(\mathrm{id}-\t)\b',\quad
		  \b'\NN=\NN\b,\quad
		  \s\b'+\b'\s=\mathrm{id},
\end{equation} 
where $\b' := \sum_{i=0}^{n-1} (-1)^i \b_i$.
The operator $\BB$ satisfies 
in general
$$\BB \BB=
(\mathrm{id}-\T)(\mathrm{id}-\t)\s\s\NN,\quad
\b\BB+\BB\b=\mathrm{id}-\T.$$ 
The 
cyclic homology 
$HC_\bullet(C)$ of a paracyclic object 
is the total homology of the bicomplex
$(E^0_{pq}:=C_{q-p}/ \mathrm{im}\,
(\mathrm{id}-\T),\b,\BB)$, 
$p,q \ge 0$ (so it only depends on the
cyclic object associated to $C$).
If $p,q$ 
take arbitrary values in 
$\mathbb{Z}$ and we consider the
direct product total complex 
we obtain the periodic 
cyclic homology $HP_\bullet(C)$.
As usual we will use the spectral sequence 
arising from filtering $E^0$ by columns as the
main tool for computing $HC_\bullet(C)$.

Recall finally that 
$\BB(D) \subset D$, where 
$D:=\mathrm{span}\{\mathrm{im}\, \s_i\}$ is the
degenerate part of $C$. Therefore,  
$\BB$ descends to the normalised complex 
$\bar C_\bullet$, where it takes the
simpler form 
$\BB=\s\NN$ since 
$\t \s=(-1)^{n+1}\t \t \s_n=-\s_0\t$.
Therefore, we
work with $\bar C_\bullet$ throughout our
explicit computations below.

\subsubsection{Twisted cyclic homology of an algebra
$\A$}\label{guruk}
If $\A$ is an algebra and $ \sigma \in \mathrm{Aut}
(\A)$, then the simplicial object $C_\bullet(\A,{}_\sigma \A)$ 
is in fact paracyclic \cite{kmt} with
$$
	\t : a_0 \otimes \ldots \otimes a_n \mapsto
	(-1)^n \sigma (a_n) \otimes a_0 
		  \otimes a_1 \otimes \ldots
		  \otimes a_{n-1}.
$$
For this paracyclic object, $\BB$ is given 
on $\bar C_\bullet(\A,{}_\sigma \A)$ by
$$
		  \BB : a_0 \otimes \ldots \otimes a_n
		  \mapsto 1 \otimes \sum_{i=0}^n (-1)^{ni}
		  \sigma (a_{n-i+1}) \otimes \ldots \otimes \sigma (a_n)
		  \otimes a_0 \otimes \ldots \otimes a_{n-i}.
$$
Following \cite{kmt}, we denote by 
$C_\bullet^\sigma(\A):=
C_\bullet(\A,{}_\sigma \A)/\mathrm{im}\,(\mathrm{id}-\T)$ 
the associated cyclic object and by
$HH^\sigma_\bullet(\A)$ and $HC^\sigma_\bullet(\A)$
its simplicial and cyclic 
homology,
respectively. By Lemma~\ref{zahn} we have
$H_\bullet(\A,{}_\sigma \A) \simeq
HH_\bullet^\sigma(\A)$
if $ \sigma $ is diagonalisable. 
This is crucial in the calculation of 
$HC^\sigma_\bullet(\A)$ since $H_\bullet(\A,{}_\sigma \A)$ is
computable via its derived functor description, while 
 $HH_\bullet^\sigma(\A)$ is the first page of the
Connes spectral sequence 
$E \Rightarrow HC^\sigma(\A)$. For an
example of an algebra and a non-diagonalisable 
$\sigma$ for which the map is not an isomorphism, see \cite{braided}, Example 3.10.

In the case $ \sigma = \mathrm{id} $,
$HC^\sigma_\bullet(\A)$ reduces to 
standard cyclic homology $HC_\bullet(\A)$
\cite{blabla,mutig}. If 
$\A=\k[X]$ for a smooth affine variety, then
the
Hochschild-Kostant-Rosenberg isomorphism identifies
$\BB$
with Cartan's exterior differential, and the
Connes spectral sequence stabilises at
$E^2$, giving 
\begin{equation}\label{chkr}
		  HP_n(\k[X]) \simeq
		\bigoplus_{i \ge 0} H_\mathrm{deRham}^{2i+n}(X),
\end{equation} 
where the right hand side is the even and odd 
algebraic de Rham
cohomology of $X$ with coefficients in $\k$, see 
e.g.~\cite{connes1,loday}.

We mention finally that
$HC_\bullet^\sigma(\A)$ can 
also be viewed as a special case of
Hopf-cyclic homology, see \cite{cm,hkrs}.

\subsection{Results}  
We now prove Theorem~\ref{themainresult}
by computing the 
twisted cyclic homology $HC^\sigma_\bullet (\A)$ of  
$\A=\k_q[SL(2)]$ for
$\sigma = \sigma _{{q^{-N},1}}$, $N \ge
2$. To do so we show that  
the spectral sequence
$E \Rightarrow HC^\sigma (\A)$
stabilises at the second page, so the
result can be read off from $E^2$. 
For the computation for all other $\sigma_{\lambda,\mu}$,
see \cite{kt}.

Lemma~\ref{zahn}
gives
$HH^\sigma_\bullet (\A) \simeq H_\bullet (\A,{}_\sigma \A)$, which is
reflected by the fact that
all the generators of $H_\bullet (\A,{}_\sigma \A)$ listed in
Section~\ref{hoho2} are invariant under 
$\T=\sigma \otimes \ldots \otimes \sigma$.
From now on we
suppress the distinction between 
$HH^\sigma_\bullet (\A)$ and $H_\bullet (\A,{}_\sigma \A)$.
We will compute  
$\BB : HH^\sigma_n(\A) \rightarrow
HH^\sigma_{n+1}(\A)$ on the vector space bases 
from
Section~\ref{hoho2}. It will then be possible
to read off directly the (co)homology
$E^2$, and it is then immediate that
$E^\bullet$ stabilises at $E^2$. 
As before, we work
throughout in the normalised complex
$\bar C_\bullet(\A,{}_\sigma \A)$. 

\subsubsection{$E^2_{pp}$ and
$HH^\sigma_1(\A) /\mathrm{im}\, \BB$} 
For $n=0$, the action of 
$\BB$ on the basis (\ref{mcquillan}) 
of $HH^\sigma_0(\A)$ 
is given by 
\begin{eqnarray}\label{ferte}
&& [1] \mapsto 0,\quad
 [x^i] \mapsto i [ x^{i-1} \otimes x ], \quad x = b,c, \; i \ge 1,\nonumber\\
&& [\omega_{N,i}] \mapsto 
		  i[\omega_{N-1,i-1} \otimes b]+
		  (N-i) [\omega_{N-1,i} \otimes c], \quad 0 \le i \le N,
		  \nonumber 
\end{eqnarray} 
(see \cite{kt}, Lemma 2.2).
Comparing 
this with the basis 
(\ref{heinsbasis}) of $HH^\sigma_1(\A)$ gives:
\begin{lemma}\label{hang}
For $p>0$ we have:
\begin{eqnarray}
		  E^2_{pp}
&=&\left\{
		  \begin{array}{ll}
		  \k[1] \quad & :\,0 \in \Slambda,\\
			0 \quad & :\,0 \notin \Slambda,
		  \end{array}\right.  
\nonumber\\ 
		  HH^\sigma_1(\A)/\mathrm{im}\, \BB
&=&\bigoplus_{i=0}^{N-2}
		  \k[\omega_{N-1,i} \otimes b] \oplus 
		  \left\{
		  \begin{array}{ll}
		  \k[b \otimes c] \quad & :\,0 \in \Slambda,\\
			0 \quad & :\,0 \notin \Slambda,
		  \end{array}\right. \nonumber
\end{eqnarray} 
where we identify elements
of $HH^\sigma_1(\A)$ with their classes
in $HH^\sigma_1(\A)/\mathrm{im}\, \BB$.
\end{lemma}

Recall that $\Slambda$ was defined in (\ref{brot}), and $0 \notin \Slambda$ if and only if $N \ge 2$ and $N$ is even. 
  
\subsubsection{$E^2_{pp+1}$ and
$HH^\sigma_2(\A) /\mathrm{im}\, \BB$} 
The basis (\ref{heinsbasis}) 
is mapped by $\BB$ to
\begin{eqnarray}
&& [b^{j-1} \otimes b],\;[c^{j-1} \otimes c]\mapsto 0,\quad
		  j>0,\>j \in \Slambda,\nonumber\\ 
&& [b \otimes c] \mapsto [1 \otimes (b \wedge c)],\quad
		  0 \in \Slambda,\nonumber\\ 
&& [\omega_{N-1,i} \otimes b] \mapsto 
		  -(N-1-i) [\omega_2^{'} (N-2,i)],\quad
		  0 \le i \le N-2,\nonumber\\ 
&& [\omega_{N-1,i} \otimes c] \mapsto
		  i [\omega_2^{'} (N-2,i-1)],\quad
		  1 \le i \le N-1.\nonumber 
\end{eqnarray} 
Now note that Lemmas ~\ref{laetst} and ~\ref{tester1} 
 imply
$
		  [1 \otimes (b \wedge c)]=0
$ 
for $0 \in \Slambda$, because in this case
$$[1 \otimes (b \wedge c)] \smallfrown 
[\partial]=0, \quad
[1 \otimes (b \wedge c)] \smallfrown 
		  [\partial']=-[b \otimes c]-[c \otimes b]=0$$
Comparing with our descriptions of 
$HH_1^\sigma(\A)/\mathrm{im}\, \BB$ 
from Lemma~\ref{hang}
and of $HH_2^\sigma (\A)$ given in 
Section~\ref{nzwei}, this yields: 
\begin{lemma}
For $p\ge1$ we have:
\begin{eqnarray}
		  E^2_{pp+1}
&=&\left\{
		  \begin{array}{ll}
		  \k[b \otimes c] \quad & :\,0 \in \Slambda,\\
			0 \quad & :\,0 \notin \Slambda,
		  \end{array}\right. \nonumber\\  
		  HH_2^\sigma (\A) / \mathrm{im}\, \BB&=&
		  \bigoplus_{i=0}^{N-2} \k[\omega_2(N-2,i)].
\nonumber
\end{eqnarray} 
\end{lemma}

\subsubsection{$E^2_{pp+2}$ and
$E^2_{pp+3}=HH^\sigma_3(\A) /\mathrm{im}\, \BB$} 
This involves a lengthier computation, so we state the
result first. 
\begin{lemma} For $\BB : HH_2^\sigma (\A) \rightarrow HH_3^\sigma (\A)$ 
we have
\begin{equation}\label{potemkin}
		  \BB([\omega_2(r,i)])=
			(2i-r)[\omega_3(r,i)],\quad
			0 \le i \le r.
\end{equation} 
Therefore, we have for $p \ge 2$
\begin{eqnarray}
		  E^2_{pp+2}
&=& \left\{
		  \begin{array}{ll}
		  \k[\omega_2(2r,r)] \quad & :\,N=2r+2,\;r \ge 0,\\
		  0 \quad & :\, \otherwise,\\
		  \end{array}\right.\nonumber\\  
		  E^2_{pp+3}
&=& HH_3^\sigma (\A) / \mathrm{im}\, \BB=\left\{
		  \begin{array}{ll}
		  \k[\omega_3(2r,r)] \quad & :\,N=2r+2,\;r \ge 0,\\
		  0 \quad & :\, \otherwise.\\
		  \end{array}\right.\nonumber  
\end{eqnarray} 
\end{lemma}
\begin{pf}
We prove (\ref{potemkin}),
the second part is then immediate. We have
 \begin{eqnarray}
&& \BB(\omega_2(r,i)) \nonumber\\ 
&=& \BB(\omega_{r,i}(
	bc \otimes (a \otimes d-d \otimes a-(q-q^{-1})c
	\otimes b) -bd \otimes (a \otimes c - qc \otimes a) 
	\nonumber\\ 
&& + da \otimes (b \otimes  c-c \otimes b) 
- q^{-1} ca \otimes (b \otimes d-qd \otimes b))) \nonumber\\
&=& 1 \otimes (\omega_{r+2,i+1} \otimes a \otimes d
		  +q^{r+2}d \otimes \omega_{r+2,i+1} \otimes a
		  +a \otimes d \otimes
		  \omega_{r+2,i+1} \nonumber\\ 
&&-\omega_{r+2,i+1} \otimes d \otimes a
		  -q^{-r-2}a \otimes \omega_{r+2,i+1} \otimes d
		  -d \otimes a \otimes
		  \omega_{r+2,i+1}\nonumber\\ 
&&-(q-q^{-1})(\omega_{r+2,i+1} \otimes c \otimes b
		  +b \otimes \omega_{r+2,i+1} \otimes c
		  +c \otimes b \otimes \omega_{r+2,i+1}) \nonumber\\ 
&&-\omega_{r,i}bd \otimes a \otimes c
		  -c \otimes \omega_{r,i}bd \otimes a
		  -q^{-r-2}a \otimes c \otimes \omega_{r,i}bd
		  \nonumber\\ 
&&+q\omega_{r,i}bd \otimes c \otimes a
		  +q^{-r-1}a \otimes \omega_{r,i}bd \otimes c
		  +q^{-r-1}c \otimes a \otimes \omega_{r,i}bd
		  \nonumber\\ 
 &&+ \omega_{r,i}da \otimes b \otimes c
		  +c \otimes \omega_{r,i}da \otimes b
		  +b \otimes c \otimes \omega_{r,i}da
		  -\omega_{r,i}da \otimes c \otimes b
		  \nonumber\\
&&-b \otimes \omega_{r,i}da \otimes c
		  -c \otimes b \otimes \omega_{r,i}da
		  -q^{-1} \omega_{r,i}ca \otimes b \otimes d
		  \nonumber\\ 
&& -q^{r+1} d \otimes \omega_{r,i}ca \otimes b
		 -q^{r+1} b \otimes d \otimes
		  \omega_{r,i}ca
		  +\omega_{r,i}ca \otimes d \otimes b \nonumber\\ 
&&+b \otimes \omega_{r,i}ca \otimes d
		  +q^{r+2} d \otimes b \otimes
		  \omega_{r,i}ca \nonumber 
\end{eqnarray}
Apply $ \smallfrown \partial$, where 
$ \partial = \frac{1}{2}(\partial^+_H+\partial^-_H)$. 
This gives
 \begin{eqnarray}
&& (\BB(\omega_2(r,i))) \smallfrown \partial \nonumber\\ 
&=& q^{r+2}d \otimes \omega_{r+2,i+1} \otimes a
		  -a \otimes d \otimes \omega_{r+2,i+1} 
		  +q^{-r-2}a \otimes \omega_{r+2,i+1} \otimes d\nonumber\\ 
&& -d \otimes a \otimes \omega_{r+2,i+1}
		  -\omega_{r,i}bd \otimes a \otimes c
		  +q^{-r-2}a \otimes c \otimes \omega_{r,i}bd
		  \nonumber\\ 
&&+q\omega_{r,i}bd \otimes c \otimes a
		  -q^{-r-1}a \otimes \omega_{r,i}bd \otimes c
		  		  +q^{-1} \omega_{r,i}ca \otimes b \otimes d
				  \nonumber\\ 
&& -q^{r+1} d \otimes \omega_{r,i}ca \otimes b
		  -\omega_{r,i}ca \otimes d \otimes b 
		  +q^{r+2} d \otimes b \otimes
		  \omega_{r,i}ca. \nonumber 
\end{eqnarray}
Now apply $ \smallfrown \frac{1}{2}(\partial^+_H-\partial^-_H)$. 
This gives
 \begin{eqnarray}
&& ((\BB(\omega_2(r,i))) \smallfrown \partial)
 \smallfrown \frac{1}{2}(\partial^+_H-\partial^-_H)\nonumber\\ 
&=& (2i-r)q^{r+2}d\omega_{r+2,i+1} \otimes a
		  +(2i-r)q^{-r-2}a\omega_{r+2,i+1} \otimes d\nonumber\\ 
&&-q^{-r-2}ac \otimes \omega_{r+1,i+1}d
		  -q\omega_{r+1,i+1}dc \otimes a\nonumber\\ 
&&-(2i-r+1)q^{-r-1}a\omega_{r+1,i+1}d \otimes c
		  +q^{-1} \omega_{r+1,i}ab \otimes d\nonumber\\ 
&&-q^{r+1}(2i-r-1) d\omega_{r+1,i}a \otimes b
		  +q^{r+2} db \otimes \omega_{r+1,i}a\nonumber\\ 
&=& (2i-r-1)q^{r+2}d\omega_{r+2,i+1} \otimes a
		  +(2i-r+1)q^{-r-2}a\omega_{r+2,i+1} \otimes d\nonumber\\ 
&&-q^{-r-2}ac \otimes \omega_{r+1,i+1}d 
		  -(2i-r+1)ad\omega_{r+1,i+1} \otimes c
		  \nonumber\\ 
&&-(2i-r-1) da\omega_{r+1,i} \otimes b
		  +q^{r+2} db \otimes \omega_{r+1,i}a. \nonumber 
\end{eqnarray}
Lemma~\ref{derifs} gives
$[\partial] \smallsmile
\frac{1}{2}[\partial^+_H-\partial^-_H]=
\frac{1}{2} [\partial^-_H] \smallsmile [\partial^+_H]$,
so by subtracting
$$
		  \b((2i-r-1)\omega_{r+2,i+1} \otimes d \otimes a
		  +b \otimes \omega_{r+1,i}a \otimes d
		  +q^{-r-2}
		  ac \otimes \omega_{r+1,i+1} \otimes d)
$$
we get in homology
\begin{eqnarray}
&& \frac{1}{2} [\BB(\omega_2(r,i))] \smallfrown 
		  ([\partial^-_H] \smallsmile [\partial^+_H])\nonumber\\ 
&=& [(2i-r-1) \omega_{r+2,i+1} \otimes bc
		  +b \otimes \omega_{r+1,i}-c \otimes \omega_{r+1,i+1}
		 \nonumber\\
&& -q^{-1}bc^2 \otimes \omega_{r+1,i+1}
		  -(2i-r+1)\omega_{r+1,i+1} \otimes c		  \nonumber\\ 
&& -(2i-r+1)q\omega_{r+3,i+2} \otimes c
		  -(2i-r-1)\omega_{r+1,i} \otimes b		  \nonumber\\ 
&& -(2i-r-1)q^{-1}\omega_{r+3,i+1} \otimes b]. \nonumber 
\end{eqnarray}
Using the calculus of
differential forms over $\k[b,c]$, that is, using 
the fact that $[f \otimes b^jc^k]=
[jf b^{j-1}c^k \otimes b]
+[kf b^jc^{k-1} \otimes c]$
for $f \in \k[b,c]$, we obtain
\begin{eqnarray}
&& \frac{1}{2} [\BB(\omega_2(r,i))] \smallfrown 
		  ([\partial^-_H] \smallsmile [\partial^+_H])\nonumber\\ 
&=& -(2i-r)[\omega_{r+1,i} \otimes b	
		  +\omega_{r+1,i+1}\otimes c]
		  \nonumber\\
&&+[((2i-r-1)-q^{-1}(r-i)-(2i-r+1)q)
		  \omega_{r+3,i+2} \otimes c\nonumber\\
&&+((2i-r-1)-(3i-r)q^{-1})\omega_{r+3,i+1} \otimes b]. \nonumber 
\end{eqnarray}
Applying
\begin{eqnarray}
&&\b(\omega_{r+1,i}a \otimes b \wedge d) = (1-q^2) \omega_{r+3,i+1} \otimes b,\nonumber\\
&&\b(\omega_{r+1,i+1}a \otimes d \wedge c) =  (q-q^{-1}) \omega_{r+3,i+2} \otimes  c\nonumber
\end{eqnarray}
the above finally simplifies to 
$$
		  \frac{1}{2} [\BB(\omega_2(r,i))] \smallfrown 
		  ([\partial^-_H] \smallsmile [\partial^+_H])=
		  -(2i-r)[\omega_{r+1,i} \otimes b	
		  +\omega_{r+1,i+1}\otimes c].
$$
The claim now
follows from Lemma~\ref{tester2}.
\end{pf}

\subsubsection{Stabilisation of the spectral 
sequence}\label{stabilo}
There is no further page of the spectral sequence to be
computed - the differential on $E^2$ maps 
$E^2_{pq}$ to $E^2_{p-2\;q+1}$, and for all $p,q$
either one space or the other is zero. Hence: 
\begin{lemma}\label{E2stable}
For $ \sigma = \sigma_{q^{-N},1}$, $N \ge 2$, 
we have
$HC_n^\sigma (\A) \simeq
\bigoplus_{p+q=n} E^2_{pq}$.
\end{lemma}
\begin{pf}
 In the case 
$0 \in \Slambda$ (i.e. $N$ odd), the $E^2$ page is as follows (the lines are for
 orientation and depict the $p=0$ and 
$q=0$ axes):
$$
\xymatrix@R=5mm@C=2mm
{& \vdots \ar@{-}[d]& \vdots & \vdots &\\
0 & 0 \ar@{-}[d] & 0 & \k[b \otimes c] & \ldots \\
0 & \bigoplus_{i=0}^{N-2} \k[\omega_2(N-2,i)]
 \ar@{-}[d] & \k[b \otimes c] & \k[1]  & \ldots\\
0 & 
{\bigoplus_{i=0}^{N-2}
\k[\omega_{N-1,i} \otimes b] \atop \oplus 
\k[b \otimes c]} \ar@{-}[d]
& \k[1] & 0 & \\
0 \ar@{-}[r] & HH_0^\sigma(\A) \ar@{-}[d]
\ar@{-}[r] & 0 \ar@{-}[r]&\\
& 0 & &}
$$
Otherwise, for $0 \notin \Slambda$ (i.e. $N$ even), the $E^2$ page is:
$$
\xymatrix@R=10mm@C=2mm
{0 & 0 \ar@{-}[d] & \k[\omega_3(N-2,\half N-1)] 
 & \k[\omega_2(N-2,\half N-1)]
 & \ldots & \\
0 & \k[\omega_3(N-2,\half N-1)] \ar@{-}[d]
 & \k[\omega_2(N-2,\half N-1)]
 & 0 & \ldots\\
0 & \bigoplus_{i=0}^{N-2} 
\k[\omega_2(N-2,i)] \ar@{-}[d]& 0 & 0 & \ldots \\
0 & \bigoplus_{i=0}^{N-2}
\k[\omega_{N-1,i} \otimes b] \ar@{-}[d]
 & 0 & 0 & \\
0 \ar@{-}[r] & HH^\sigma_0(\A) \ar@{-}[d] \ar@{-}[r] & 0 \ar@{-}[r] & & \\
& 0 &  & &
}
$$
\end{pf}

Theorem \ref{themainresult} now follows
immediately.

\section{Cyclic cohomology}\label{cyccohom}
\subsection{Background}\label{woe}
In this final section we finish the
proof of Theorem \ref{themainresult2}
which we began in Section \ref{como} :
we construct a twisted cyclic
3-cocycle $\ccc$ that pairs nontrivially with $\dA$,
and hence represents a generator in Connes' $ \lambda $-complex  
of $HC_{\sigma_{q^{-2},1}}^3(\A) \simeq \k$.

\subsubsection{Cyclic cocycles}\label{rano}
In view of (\ref{skeip}),
$ \mathrm{im}(\mathrm{id} - \t) \subset C$ is for any cyclic
$\k$-vector space a subcomplex with
respect to $\b$. As Connes 
showed, cyclic homology
can be realised  
as the homology of the quotient:
$$
		  HC_\bullet(C) \simeq H_\bullet(C/\mathrm{im}(\mathrm{id}-\t),\b).
$$
For $C=C^\sigma(\A)$ we can dually consider 
Hochschild cochains
$ \varphi \in C^n(\A,({}_\sigma \A)^*)$ which are
twisted cyclic, that is, which satisfy 
\begin{equation}\label{cyclisch}
		  \varphi(a_0,\ldots,a_n)=
		  (-1)^n \varphi(\sigma(a_n),a_0,\ldots,a_{n-1})
\end{equation} 
for all 
$a_0 ,\ldots,a_n \in \A$. These form a subcomplex
of $(C^\bullet(\A,({}_\sigma \A)^*),\b)$ whose
cohomology is twisted cyclic cohomology 
$HC_\sigma^\bullet(\A) \simeq
(HC^\sigma_\bullet(\A))^*$.

If we work with the normalised complex, then
$ \varphi (a_0,\ldots,a_n)=0$ whenever $a_i \in \k$ for
some $i>0$. For cyclic $ \varphi $ this property
obviously extends to $i=0$. Conversely, a Hochschild
cocycle that vanishes on 
$1 \otimes a_1 \otimes \ldots
\otimes a_n$ is twisted cyclic as follows by applying it 
to $\b(1 \otimes a_0 \otimes \ldots \otimes a_n)$. 

\subsection{Results} 
The natural question is whether
the Hochschild 3-cocycle 
$\varphi$ defined in (\ref{losfundamentos}) 
 is already cyclic. We show that it is not, then construct a coboundary 
$ \eta$ by which $\varphi$ differs from
a cyclic cocycle.

\subsubsection{$\varphi$ is not cyclic} 
 As remarked at the end of 
Section~\ref{rano}, cyclicity  is equivalent to the condition
that 
\begin{equation}\label{cyclicornot}
		  \varphi(1,a_1,a_2,a_3)=0
\end{equation}
for all $a_1,a_2,a_3 \in \A$.
 Now, for all $i,j \ge 0$ we have 
\begin{eqnarray}
&& \sigma_{q,q{-1}}( \partial_E^+(d^jc))\;
		  \partial^+_F(a^ib)=q^{i-2j}
		  d^{j+1}a^{i+1},\nonumber\\ 
&& \sigma_{q^2,q^{-2}}
		  (\partial^+_H(e_{j-i,0,0}))=
		  (i-j)q^{2(j-i)}e_{j-i,0,0},\nonumber
\end{eqnarray} 
and therefore
$$
		  \varphi(1,e_{j-i,0,0},d^jc,a^ib)=
		  q^{-i}(i-j) \int_{[1]} e_{j-i,0,0}\;
		  d^{j+1}a^{i+1}=q^{-i}(i-j).
$$
Hence (\ref{cyclicornot}) fails and so $\varphi$ is not cyclic.

\subsubsection{The correction term}\label{letzz} 
We make the ansatz
$$
		  \eta(\cdot):=\int_?
		  \cdot \smallfrown
		  (\partial^+_H \smallsmile 
		  (\sigma_{\lambda_1,\mu_1}-\mathrm{id})\smallsmile 
		  (\sigma_{\lambda_2,\mu_2}-\mathrm{id})),
$$
where $\int_?$ is a suitable twisted
trace and $\lambda_1,\mu_1,\lambda_2,\mu_2$
are complex parameters.
\begin{lemma}
If 
$\lambda_1=\lambda_2=1,\mu_2=\mu_1^{-1}
 \neq 1$ and
$$
		  \int_?=-\frac{\mu_2}{(\mu_2-1)^2}\int_{[bc]} \in
		  H^0(\A,({}_{\sigma_{q^{-2},1}} \A)^*),
$$
then $ \varphi +\eta$ is cyclic 
and as a Hochschild cocycle
cohomologous to  
$ \varphi$.
\end{lemma}
\begin{pf}
For any automorphism $ \sigma $ of an algebra, 
$ \sigma - \mathrm{id}$ is a $ \sigma $-twisted derivation which
 is inner (it is simply the twisted commutator 
with $1 \in \A$). Therefore its cohomology class
 vanishes, 
 so $ \eta$ is
 automatically a coboundary.
Furthermore,  for all $i,j \ge 0$ we have
\begin{eqnarray}
&& \sigma_{\lambda_2,\mu_2}
		  ((\sigma_{\lambda_1,\mu_1}-\mathrm{id})(d^jc))
		  (\sigma_{\lambda_2,\mu_2}-\mathrm{id})(a^ib)\nonumber\\ 
&=& q^{-i}\lambda_2^{-j}\mu_2^{-1}
		  (\lambda_1^{-j}\mu_1^{-1}-1)
		  (\lambda_2^i \mu_2-1) d^ja^ibc,\nonumber\\ 
&& \sigma_{\lambda_1 \lambda_2,\mu_1 \mu_2}
		  (\partial^+_H(e_{j-i,0,0}))\nonumber\\ 
&=&  (i-j)(\lambda_1 \lambda_2)^{j-i}e_{j-i,0,0}.\nonumber
\end{eqnarray} 
Hence
\begin{eqnarray}
&& \eta(1,e_{j-i,0,0},d^jc,a^ib) \nonumber\\ 
&=& q^{-i}\lambda_2^{-j}\mu_2^{-1}
		  (\lambda_1^{-j}\mu_1^{-1}-1)
		  (\lambda_2^i \mu_2-1)
		  (i-j)(\lambda_1 \lambda_2)^{j-i}
		  \int_?  e_{j-i,0,0}\;d^ja^ibc.\nonumber 
\end{eqnarray}
Therefore, we need $\int_?$ to be a
$\sigma_{q^{-2} \lambda_1 \lambda_2,\mu_1 \mu_2}$-twisted
trace for which
$$
		  \int_?  e_{j-i,0,0}\;d^ja^ibc = 
		  - \frac{\lambda_1^{i-j} \lambda_2^i
		  \mu_2}{(\lambda_1^{-j}\mu_1^{-1}-1)
		  (\lambda_2^i \mu_2-1)
		  } .
$$
It is easily checked that  
$ \varphi(1,e_{i,j,k},e_{l,m,n},e_{r,s,t})=
\eta(1,e_{i,j,k},e_{l,m,n},e_{r,s,t})=0$
for all other $i,j,k,l,m,n,r,s,t$, 
since $b$ and $c$ are twisted central.
\end{pf}

This completes the proof of Theorem \ref{themainresult2}.

\end{document}